\documentclass[10pt]{confm}
\usepackage{inputenc}
\inputencoding{cp1251}
\usepackage[russian,english]{babel}
\usepackage{nauka}

\newcounter{theorem}
\newcommand{\theor}{\vspace{3mm}\noindent\refstepcounter{theorem}%
{\bf Theorem \thetheorem.}\,\,}
\newcommand{\lr}[1]{\left(#1\right)}

\begin{document}

\setcounter{page}{1} \thispagestyle{empty} \ab

%\begin{titlepage}

\title{BUSCHMAN--ERDELYI \,TRANSMUTATIONS,\protect\\[2.0mm]
CLASSIFICATION AND  \,APPLICATIONS}

\author{S.M.~Sitnik}
\bigskip

\begin{center}
This work is dedicated to Professor  Sergei V. Rogosin on the remarkable occasion of his 60th birthday.
\end{center}
\bigskip

\titlehead{TRANSMUTATIONS \,AND \,APPLICATIONS}
\authorhead{S.M.~Sitnik}

\address{Voronezh Institute of the Ministry of Internal Affairs of Russia,
\protect\\ 53, Patriotov prospekt,  394065, Voronezh,
Russia}

%\tableofcontents

%%\received{}
%%e-mail: box2008in@gmail.com

{\small\baselineskip 10pt  \abstract{
\ \ This survey paper contains brief historical information, main known facts and original author's results on the theory of  transmutations and some applications. Operators of  Buschman--Erdelyi type were first studied by E.T.~Copson, R.G.~Buschman and A.~Erdelyi as integral operators. In 1990's the author was first to prove transmutational nature of these operators and published papers with detailed  study of their properties. This class include as special cases such famous objects as Sonine--Poisson--Delsarte transmutations and fractional Riemann--Lioville integrals. In this paper Buschman--Erdelyi transmutations are fully classified as operators of the first kind with special case of zero order smoothness operators, second kind and third kind  with special case of unitary Sonine--Katrakhov and Poisson--Katrakhov transmutations. We study such properties as transmutational conditions, factorizations, norm estimates, connections with classical integral transforms. Applications are considered to singular partial differential equations, embedding theorems with sharp constants in Kipriyanov spaces, Euler--Poisson--Darboux equation including Copson lemma, generalized translations, Dunkl operators,  Radon transform, generalized harmonics theory, Hardy operators, V.\,Katrakhov's results on pseudodifferential operators  and problems  of new kind  for equations with solutions of arbitrary growth at singularity.}


}

\newpage

\addcontentsline{toc}{section}{-------------------------------------------------TABLE OF CONTENTS---------------------------------------------}

\section{
\\Introduction: an idea of transmutations, historical information and applications.}

\subsection{Transmutation operators.}

\medskip

Transmutation theory is an essential generalization of matrix similarity theory. Let start with the main definition.

Definition 1. For a given pair of operators $(A,B)$ an operator $T$ is called transmutation (or intertwining) operator if on  elements of some functional spaces the next property is valid
\begin{equation}
\label{1.1}
\Large {T\,A=B\,T.}
\end{equation}

It is obvious that the notion of transmutation is direct and far going generalization of similarity notion from linear algebra. But  transmutations do not reduce to similar operators because intertwining operators often are not bounded in classical spaces and the inverse operator may be not exist or bounded in the same space. As a consequence spectra of   intertwining operators are not the same as a rule. Moreover transmutations may be unbounded. It is the case for Darboux transformations which are defined for a pair of differential operators and are differential operators themselves, in this case all three operators are unbounded in classical spaces. But the theory of Darboux transformations is included in transmutation theory too. Also a pair of intertwining operators may not be differential ones. In transmutation theory there are problems for next varied types of operators: integral, integro--differential, difference--differential (e.g.  the Dunkl operator), differential or integro--differential of infinite order (e.g. in connection with Schur's lemma), general linear operators in functional spaces, pseudodifferential and abstract  differential operators.

All classical integral transforms due to the definition 1 are also special cases of transmutations, they include Fourier, Petzval (Laplace), Mellin, Hankel, Weierstrass, Kontorovich--Lebedev, Meyer, Stankovic, finite Grinberg and other transforms.

In quantum physics in study of Shr\"{o}dinger equation and inverse scattering theory underlying transmutations are called wave operators.

Commuting operators are also a special class of transmutations. The most important class consists of operators commuting with derivatives. In this case transmutations as commutants are usually in the form of formal series, pseudodifferential or infinite order differential operators. Finding of commutants is directly connected with finding all  transmutations in the given functional space. For these problems works a  theory of operator convolutions, including Berg--Dimovski convolutions. Also more and more applications are developed connected with transmutation theory for commuting differential operators, such problems are based on classical results of J.L.\,Burchnall, T.W.\,Chaundy. Transmutations are also connected with factorization problems for integral and differential operators. Special class of transmutations are the so called Dirichlet--to--Neumann and  Neumann--to--Dirichlet operators which link together solutions of the same equation but with different kinds of boundary conditions.

And how transmutations usually works? Suppose we study properties for a  rather complicated operator
$A$. But suppose also that we know corresponding  properties for a  model more simple operator $B$ and transmutation (\ref{1.1}) readily exists. Then we usually may copy results for the model operator $B$ to a more complicated operator $A$. This is shortly the main idea of transmutations.

Let us for example consider an equation $Au=f$, then applying to it a transmutation with property (\ref{1.1}) we consider a new equation $Bv=g$, with $v=Tu, g=Tf$. So if we can  solve the simpler   equation $Bv=g$ then the initial one is also solved and has solution $u=T^{-1}v$. Of course it is supposed that the  inverse operator exist and its explicit form is known. This is a simple application of transmutation technique for proving formulas for solutions of ordinary and partial differential equations.

The next monographs  \cite{Car1}-\cite{Tri1} are completely devoted to the transmutation theory and its applications. Moreover essential parts of  monographs \cite{CarSho}-\cite{Hro} include material on transmutations, the complete list of  books which consider some transmutational problems is now near of 100 items.

We specially distinguish the book \cite{FaNa}. In it the most difficult problems of transmutation theory were solved. Among them an existence of transmutations was proved for high order differential equations
with variable coefficients including correcting errors of previous papers of Delsarte and Lions, the complete theory of Bianchi equation, extension of V.\,Marchenko theory of operator--analytic functions, results on operators commuting in spaces of analytic functions.

We use the term "transmutation"\  due to \cite{Car3}: "Such operators are often called transformation operators by the Russian school (Levitan, Naimark, Marchenko et. al.), but transformation seems too broad a term, and, since some of the machinery seems "magical"\  at times, we have followed Lions and Delsarte in using the word "transmutation".

Now the transmutation theory is a completely formed part of mathematical world in which methods and ideas from different areas are used: differential and integral equations, functional analysis, function theory, complex analysis, special functions, fractional integrodifferentiation.

Transmutation theory is deeply connected with many applications in different fields of mathematics. Transmutations are applied in inverse problems via the generalized Fourier transform, spectral function and famous Levitan equation; in scattering theory  the Marchenko equation is formulated in terms of transmutations; in spectral theory transmutations help to prove trace formulas and asymptotics for spectral function; estimates for transmutational kernels control stability in inverse and scattering problems;
for nonlinear equations via Lax method transmutations for Sturm--Lioville problems lead to proving existence and explicit formulas for soliton solutions. Special kinds of transmutations are generalized analytic functions, generalized translations and convolutions, Darboux transformations. In the theory of partial differential equations transmutations works for proving explicit correspondence formulas among solutions of perturbed and non--perturbed equations, for singular and degenerate equations, pseudodifferential operators, problems with essential singularities at inner or corner points, estimates of solution decay for elliptic and ultraelliptic equations. In function theory transmutations are applied to embedding theorems and generalizations of Hardy operators, Paley--Wiener theory, generalizations of harmonic analysis based on generalized translations. Methods of transmutations are used in many applied problems: investigation of Jost solutions in scattering theory, inverse problems, Dirac and other matrix systems of differential equations, integral equations with special function kernels, probability theory and random processes, stochastic random equations, linear  stochastic estimation, inverse problems of geophysics and transsound gas dynamics. Also a number of applications of transmutations to nonlinear equations is permanently increased.

 In fact the modern transmutation theory originated from two basic examples \cite{Sit1}.
The first is transmutations $T$ for Sturm--Lioville problems with some potential $q(x)$ and natural boundary conditions
\begin{equation}
T(D^2\,y(x)+q(x)y(x))=D^2\,(Ty(x)), D^2\,y(x)=y''(x),
\end{equation}

The second example is a problem of studying transmutations intertwining the Bessel operator $B_\nu$ and the second derivative:
\begin{equation}
T\lr{B_\nu} f=\lr{D^2} Tf,\ \  B_{\nu}=D^2+\frac{2\nu +1}{x}D,\ \  D^2=\frac{d^2}{dx^2},\ \  \nu \in \mathbb{C}.
\end{equation}
This class of transmutations includes Sonine--Poisson--Delsarte, Buschman--Erdelyi operators and generalizations. Such transmutations found many applications for a special class of partial differential equations with singular coefficients. A typical equation of this class is  the $B$--elliptic equation with the Bessel operator in some variables of the form
\begin{equation}
\label{59}
\sum_{k=1}^{n}B_{\nu,x_k}u(x_1,\dots, x_n)=f.
\end{equation}
Analogously  $B$--hyperbolic and $B$--parabolic equations are considered, this terminology was proposed by I.\,Kipriyanov. This class of equations was first studied by Euler, Poisson, Darboux and continued in Weinstein's theory of generalized axially symmetric potential (GASPT). These problems were further investigated by Zhitomirslii, Kudryavtsev, Lizorkin, Matiychuk, Mikhailov, Olevskii, Smirnov, Tersenov, He Kan Cher, Yanushauskas, Egorov and others.

In the most detailed and complete way  equations with Bessel operators were studied by the Voronezh mathematician I.A.\,Kipriyanov and his disciples Ivanov, Ryzhkov, Katrakhov, Arhipov, Baidakov, Bogachov, Brodskii, Vinogradova, Zaitsev, Zasorin, Kagan, Katrakhova, Kipriyanova, Kononenko, Kluchantsev, Kulikov, Larin, Leizin, Lyakhov, Muravnik, Polovinkin, Sazonov, Sitnik, Shatskii, Yaroslavtseva. The essence of Kipriyanov's school results was published in  \cite{Kip1}. For classes of equations with Bessel operators I.\,Kipriyanov introduced  special functional spaces which were named after him \cite{Kip2}.
In this field interesting results were investigated by Katrakhov and his disciples, now these problems are considered by  Gadjiev, Guliev, Glushak, Lyakhov with their coauthors and students. Abstract equations of the form (\ref{59}) originated from the monograph \cite{CarSho} were considered by  Egorov, Repnikov, Kononenko, Glushak, Shmulevich and others. And transmutations are one of basic tools for equations with Bessel operators, they are applied to construction of solutions, fundamental solutions, study of singularities, new boundary--value and other problems.

Some words about the structure of this publication.
This is a survey article on transmutations of special classes. But the main result on norm estimates and unitarity of Buschman--Erdelyi transmutations is completely proved  (theorem \ref{tnorm}) as many other facts are consequences of this theorem. In the first section historical and priority information is provided. An author's classification of different classes of  Buschman--Erdelyi transmutations is introduced. Based on this classification  Buschman--Erdelyi transmutations of the first kind and zero order operators are studied in the second section.  Buschman--Erdelyi transmutations of the second kind are considered in the third section. In the fourth section  Buschman--Erdelyi transmutations of the third kind and also Sonine--Katrakhov and Poisson--Katrakhov unitary transmutations are considered. In the final fifth section different applications of  Buschman--Erdelyi transmutations are listed, mostly inevitably briefly.
They include embedding theorems for Kipriyanov spaces, solution representations for partial differential equations with Bessel operators, Euler--Poisson--Darboux equations and Copson's lemma for them, generalized translations, Dunkl operators, Radon transform, generalized spherical harmonics and $B$--harmonic polynomials, unitarity for some generalizations of Hardy operators. In the final part of this section some results of V.\,Katrakhov is mentioned on a new class of pseudodifferential operators and remarkable problems introduced by him with   $K$---trace for solutions with infinite order singularities.

Also we must note that the term "operator"\ is used in this paper for brevity in the broad and sometimes not exact meaning, so appropriate domains and function classes are not always specified. It is easy to complete and make strict for every special result.

\medskip
\subsection{Buschman--Erdelyi transmutations.}
\medskip

The term "Buschman--Erdelyi transmutations"\ was introduced by the author and is now accepted. Integral
equations with these operators were studied in mid--1950th. The author was first to prove the transmutational nature of these operators. The classical Sonine and Poisson operators are special cases of Buschman--Erdelyi transmutations and Sonine--Dimovski and Poisson--Dimovski transmutations are their generalizations for hyper--Bessel equations and functions.

Buschman--Erdelyi transmutations have many modifications. The author introduced convenient classification of them. Due to this classification we introduce Buschman--Erdelyi transmutations of the first kind, their kernels are expressed in terms of Legendre functions of the first kind. In the limiting case we define Buschman--Erdelyi transmutations of zero order smoothness being important in applications. Kernels of  Buschman--Erdelyi transmutations of the second  kind are expressed in terms of Legendre functions of the second kind. Some combination of operators of  the first kind and the second kind leads to  operators of the third kind. For the special choice of parameters they are unitary operators in the standard Lebesgue space.  The author proposed terms
"Sonine--Katrakhov"\ and "Poisson--Katrakhov"\ transmutations in honor of V.\,Katrakhov who  introduced and studied  these operators.

The study of integral equations and invertibility for Buschman--Erdelyi operators  was started in 1960-th by P.Buschman and A.Erdelyi  \cite{Bu1}--\cite{Er2}. These operators also were investigated by Higgins, Ta Li, Love, Habibullah, K.N.\,Srivastava, Ding Hoang An, Smirnov, Virchenko, Fedotova, Kilbas, Skoromnik and others.
During this period for this class of operators were considered only problems of solving integral equations, factorization and invertibility, cf.  \cite{KK}.

The most detailed study of Buschman--Erdelyi transmutations was taken by the author in 1980--1990th
\cite{Sit3}--\cite{Sit4} and continued in \cite{Sit2}--\cite{Sit10} and some other papers.
Interesting results were proved by N.Virchenko and A.Kilbas and their disciples \cite{KiSk1}--\cite{KiSk2}, \cite{ViFe}.

\section{\\Buschman--Erdelyi transmutations of the first kind.}

\setcounter{equation}{3}

\subsection{Sonine--Poisson--Delsarte transmutations.}

Let us first consider  the most well--known transmutations for the Bessel operator and the second derivative:
\begin{equation}
\label{2.1}
T\lr{B_\nu} f=\lr{D^2} Tf, B_{\nu}=D^2+\frac{2\nu +1}{x}D, D^2=\frac{d^2}{dx^2}, \nu \in \mathbb{C}.
\end{equation}

Definition 2.
The Poisson transmutation is defined by
\begin{equation}
\label{2.2}
P_{\nu}f=\frac{1}{\Gamma(\nu+1)2^{\nu}x^{2\nu}}
\int_0^x \left( x^2-t^2\right)^{\nu-\frac{1}{2}}f(t)\,dt,\Re \nu> -\frac{1}{2}.
\end{equation}
The Sonine transmutation is defined by
\begin{equation}
\label{2.3}
S_{\nu}f=\frac{2^{\nu+\frac{1}{2}}}{\Gamma(\frac{1}{2}-\nu)}\frac{d}{dx}
\int_0^x \left( x^2-t^2\right)^{-\nu-\frac{1}{2}}t^{2\nu+1}f(t)\,dt,\Re \nu< \frac{1}{2}.
\end{equation}
Operators (\ref{2.2})--(\ref{2.3}) intertwine by formulas
\begin{equation}
\label{56}
S_\nu B_\nu=D^2 S_\nu,\ P_\nu D^2=B_\nu P_\nu.
\end{equation}
The definition may be extended to $\nu\in\mathbb{C}$.
We will use more historically exact term Sonine--Poisson--Delsarte transmutations \cite{Lev3}.

An important generalization for Sonine--Poisson--Delsarte are transmutations for hyper--Bessel functions.
Such functions were first considered by Kummer and Delerue. The detailed study was done by Dimovski and his coauthors \cite{Dim}. These transmutations are called Sonine--Dimovski and  Poisson--Dimovski by Kiryakova    \cite{Kir}. In hyper--Bessel functions theory the leading role is for Obreshkoff integral transform \cite{Kir}. It is a transform with Mayer's $G$--function kernel which generalize Laplace, Mellin, sine and cosine Fourier, Hankel,  Mayer and other classical transforms. Different results on hyper--Bessel functions, connected equations and transformed were many times reopened. The same is true for the Obreshkoff integral transform. It my opinion the Obreshkoff  transform together  with Fourier, Mellin, Laplace, Stankovic transforms are basic elements from which many other transforms are constructed with corresponding applications.

\medskip
\subsection{Definition and main properties of Buschman--Erdelyi transmutations of the first kind.}
\medskip

Let define and study main properties of  Buschman--Erdelyi transmutations of the first kind.
This class of transmutations for some choice of parameters generalizes Sonine--Poisson--Delsart  transmutations, Riemann--Liouville  and Erdelyi--Kober fractional integrals, Mehler--Fock transform.

Definition 3.
Define Buschman--Erdelyi operators of the first kind by
\begin{eqnarray}
\label{71}
B_{0+}^{\nu,\mu}f=\int_0^x \left( x^2-t^2\right)^{-\frac{\mu}{2}}P_\nu^\mu \left(\frac{x}{t}\right)f(t)d\,t,\\
E_{0+}^{\nu,\mu}f=\int_0^x \left( x^2-t^2\right)^{-\frac{\mu}{2}}\mathbb{P}_\nu^\mu \left(\frac{t}{x}\right)f(t)d\,t,\\
B_{-}^{\nu,\mu}f=\int_x^\infty \left( t^2-x^2\right)^{-\frac{\mu}{2}}P_\nu^\mu \left(\frac{t}{x}\right)f(t)d\,t,\\
\label{72}
E_{-}^{\nu,\mu}f=\int_x^\infty \left( t^2-x^2\right)^{-\frac{\mu}{2}}\mathbb{P}_\nu^\mu \left(\frac{x}{t}\right)f(t)d\,t.\\ \nonumber
\end{eqnarray}
here $P_\nu^\mu(z)$ is the Legendre function of the first kind , $\mathbb{P}_\nu^\mu(z)$ is this function on the cut $-1\leq  t \leq 1$, $f(x)$  is a locally summable function with some growth conditions at $x\to 0,x\to\infty$. Parameters $\mu,\nu\in\mathbb{C}$, $\Re \mu <1$,  $\Re \nu \geq -1/2$.

Now consider main properties  for this class of transmutations following essentially \cite{Sit3}, \cite{Sit4}, and also \cite{Sit1}, \cite{Sit2}. All functions further are defined on positive semiaxis.
So we use notations $L_2$ for the functional space $L_2(0, \infty)$ and $L_{2, k}$ for power weighted
space $L_{2, k}(0, \infty)$ equipped with norm
\begin{equation}
\int_0^\infty |f(x)|^2 x^{2k+1}\,dx.
\end{equation}
$\mathbb{N}$ \ denote set of natural, $\mathbb{N}_0$--positive integer, $\mathbb{Z}$--integer and  $\mathbb{R}$--real numbers.

First add to definition 3 a case of parameter $\mu =1$.

Definition 4.
Define for $\mu =1$ Buschman--Erdelyi operators of zero order smoothness by
\begin{eqnarray}
\label{73}
B_{0+}^{\nu,1}f=\frac{d}{dx}\int_0^x P_\nu \left(\frac{x}{t}\right)f(t)\,dt,\\
\label{731}
E_{0+}^{\nu,1}f=\int_0^x P_\nu \left(\frac{t}{x}\right)\frac{df(t)}{dt}\,dt,\\
\label{732}
B_{-}^{\nu,1}f=\int_x^\infty P_\nu \left(\frac{t}{x}\right)(-\frac{df(t)}{dt})\,dt,\\
\label{733}
E_{-}^{\nu,1}f=(-\frac{d}{dx})\int_x^\infty P_\nu \left(\frac{x}{t}\right)f(t)\,dt,
\end{eqnarray}
here $P_\nu(z)=P_\nu^0(z)$ is the Legendre function.

\theor{\it
The next formulas hold true for factorizations of Buschman--Erdelyi transmutations for suitable functions via Riemann--Liouville fractional integrals and  Buschman--Erdelyi operators of zero order smoothness:

 \begin{equation}\label{1.9}
{B_{0+}^{\nu,\,\mu} f=I_{0+}^{1-\mu}~ {_1 S^{\nu}_{0+}f},~B_{-}^{\nu, \,\mu} f={_1 P^{\nu}_{-}}~ I_{-}^{1-\mu}f,}
\end{equation}
\begin{equation}\label{1.10}
{E_{0+}^{\nu,\,\mu} f={_1 P^{\nu}_{0+}}~I_{0+}^{1-\mu}f,~E_{-}^{\nu, \, \mu} f= I_{-}^{1-\mu}~{_1 S^{\nu}_{-}}f.}
\end{equation}
}

These formulas allow to separate parameters $\nu$ and $\mu$. We will prove soon that operators
\eqref{73}--\eqref{733} are isomorphisms of $L_2(0, \infty)$ except for some special parameters. So operators  \eqref{71}--\eqref{72} roughly speaking are of the same smoothness in $L_2$ as integrodifferentiations $I^{1-\mu}$ and they coincide with them for $\nu=0$. It is also possible to define Buschman--Erdelyi operators for all $\mu\in\mathbb{C}$.

Definition 5.{ Define the number $\rho=1-Re\,\mu $ as smoothness order for Buschman--Erdelyi operators \eqref{71}--\eqref{72}.}

So for $\rho > 0$ (otherwise for $Re\, \mu > 1$) Buschman--Erdelyi operators are smoothing and for $\rho < 0$ (otherwise for $Re\, \mu < 1$) they decrease smoothness in $L_2$ spaces. Operators \eqref{73}--\eqref{733} for which $\rho = 0$ due to definition 5 are of zero smoothness order.

For some special parameters $\nu,~\mu$ Buschman--Erdelyi operators of the first kind are reduced to other known operators. So for $\mu=-\nu$ or $\mu=\nu+2$ they reduce to Erdelyi--Kober operators, for $\nu = 0$ they reduce to fractional integrodifferentiation $I_{0+}^{1-\mu}$ or $I_{-}^{1-\mu}$, for  $\nu=-\frac{1}{2}$, $\mu=0$ or $\mu=1$ kernels reduce to elliptic integrals, for $\mu=0$,  $x=1$, $v=it-\frac{1}{2}$ the operator $B_{-}^{\nu, \, 0}$ differs only by a constant from Mehler--Fock transform.

As a pair for the Bessel operator consider a connected one
\begin{equation}
\label{75}
L_{\nu}=D^2-\frac{\nu(\nu+1)}{x^2}=\left(\frac{d}{dx}-\frac{\nu}{x}\right)
\left(\frac{d}{dx}+\frac{\nu}{x}\right),
\end{equation}
which for $\nu \in \mathbb{N}$ is an angular momentum operator from quantum physics. Their transmutational relations are established in the next theorem.

\theor \label{tOP}
{\it For a given pair of transmutations $X_\nu, Y_\nu$
\begin{equation}
\label{76}
X_\nu L_{\nu}=D^2 X_\nu , Y_\nu D^2 = L_{\nu} Y_\nu
\end{equation}
define the new pair of transmutations by formulas
\begin{equation}
\label{77}
S_\nu=X_{\nu-1/2} x^{\nu+1/2}, P_\nu=x^{-(\nu+1/2)} Y_{\nu-1/2}.
\end{equation}
Then for the new pair  $S_\nu, P_\nu$ the next formulas are valid:
\begin{equation}
\label{78}
S_\nu B_\nu = D^2 S_\nu, P_\nu D^2 = B_\nu P_\nu.
\end{equation}
}

\theor\label{t6}
{\it Let $Re \, \mu \leq 1$. Then an operator $B_{0+}^{\nu, \, \mu}$ is a  Sonine type transmutation and
\eqref{76} is valid.
}

The same result holds true for other Buschman--Erdelyi operators, $E_{-}^{\nu, \, \mu}$ is Sonine type
and $E_{0+}^{\nu, \, \mu}$,  $B_{-}^{\nu, \, \mu}$ are Poisson type transmutations.

From these transmutational connections we conclude that  Buschman--Erdelyi operators link corresponding eigenfunctions for two operators. They lead to formulas for Bessel functions via exponents and trigonometric functions and vice versa which generalize classical Sonine and Poisson formulas.

Now consider factorizations of Buschman--Erdelyi operators.
First let list main forms of fractional integrodifferentiations: Riemann--Liouville, Erdelyi--Kober, fractional integral by function $g(x)$, cf. \cite{KK}.
\begin{eqnarray}
\label{61}
I_{0+,x}^{\alpha}f=\frac{1}{\Gamma(\alpha)}\int_0^x \left( x-t\right)^{\alpha-1}f(t)d\,t,\\ \nonumber
I_{-,x}^{\alpha}f=\frac{1}{\Gamma(\alpha)}\int_x^\infty \left( t-x\right)^{\alpha-1}f(t)d\,t,
\end{eqnarray}
\begin{eqnarray}
\label{62}
I_{0+,2,\eta}^{\alpha}f=\frac{2 x^{-2\lr{\alpha+\eta}}}{\Gamma(\alpha)}\int_0^x \left( x^2-t^2\right)^{\alpha-1}t^{2\eta+1}f(t)d\,t,\\ \nonumber
I_{-,2,\eta}^{\alpha}f=\frac{2 x^{2\eta}}{\Gamma(\alpha)}\int_x^\infty \left( t^2-x^2\right)^{\alpha-1}t^{1-2\lr{\alpha+\eta}}f(t)d\,t,
\end{eqnarray}
\begin{eqnarray}
\label{63}
I_{0+,g}^{\alpha}f=\frac{1}{\Gamma(\alpha)}\int_0^x \left( g(x)-g(t)\right)^{\alpha-1}g'(t)f(t)d\,t,\\ \nonumber
I_{-,g}^{\alpha}f=\frac{1}{\Gamma(\alpha)}\int_x^\infty \left( g(t)-g(x)\right)^{\alpha-1}g'(t)f(t)d\,t,
\end{eqnarray}
in all cases $\Re\alpha>0$ and operators may be further defined for all $\alpha$ \cite{KK}.
In case of $g(x)=x$  (\ref{63}) reduces to Riemann--Liouville, in case of $g(x)=x^2$ (\ref{63}) reduces to Erdelyi--Kober and in case of $g(x)=\ln x$ to Hadamard fractional integrals.

\theor\label{tfact7} The next factorization formulas are valid for Buschman--Erdelyi operators of the first kind via Riemann--Liouville and Erdelyi--Kober fractional integrals
\begin{eqnarray}
& & B_{0+}^{\nu, \, \mu}=I_{0+}^{\nu+1-\mu} I_{0+; \, 2, \, \nu+ \frac{1}{2}}^{-(\nu+1)} {\lr{\frac{2}{x}}}^{\nu+1}\label{2.17}, \\
& & E_{0+}^{\nu, \, \mu}= {\lr{\frac{x}{2}}}^{\nu+1} I_{0+; \, 2, \, - \frac{1}{2}}^{\nu+1} I_{0+}^{-(\nu+\mu)}  \label{2.18}, \\
& & B_{-}^{\nu, \, \mu}= {\lr{\frac{2}{x}}}^{\nu+1}I_{-; \, 2, \, \nu+ 1}^{-(\nu+1)} I_{-}^{\nu - \mu+2}  \label{2.19}, \\
& & E_{-}^{\nu, \, \mu}= I_{-}^{-(\nu+\mu)} I_{-; \, 2, \, 0} ^{\nu+1} {\lr{\frac{x}{2}}}^{\nu+1}  \label{2.20}.
\end{eqnarray}

Sonine--Poisson--Delsarte transmutations also are special cases for this class of operators.

Now let study properties of Buschman--Erdelyi operators of zero order smoothness defined by (\ref{73}).
A similar operator was introduced by Katrakhov by multiplying the Sonine operator with fractional integral, his aim was to work with transmutation obeying good estimates in $L_2(0,\infty)$.

We use the Mellin transform defined by \cite{Mar}
\begin{equation}
\label{710}
g(s)=M{f}(s)=\int_0^\infty x^{s-1} f(x)\,dx.
\end{equation}
The Mellin convolution is defined by
\begin{equation}
\label{711}
(f_1*f_2)(x)=\int_0^\infty  f_1\left(\frac{x}{y}\right) f_2(y)\,\frac{dy}{y},
\end{equation}
so the convolution operator with kernel $K$ acts under Mellin transform as multiplication on multiplicator \begin{eqnarray}
\label{con}
M[Af](s)=M\ [\int_0^\infty  K\left(\frac{x}{y}\right) f(y)\,\frac{dy}{y}](s)=M[K*f](s)
=m_A(s)M{f}(s),\\\nonumber m_A(s)=M[K](s).\phantom{1111111111111111111}
\end{eqnarray}

We observe that Mellin transform is a generalized Fourier transform on semiaxis with Haar measure
$\frac{dy}{y}$ \cite{Hel}. It plays important role for special functions, for example the gamma function is a Mellin transform of the exponential. With Mellin transform the important
breakthrough in evaluating integrals was done in 1970th when mainly by O.Marichev the famous Slater's theorem was adapted for calculations. The Slater's theorem taking the Mellin transform as input give the function itself as output via hypergeometric functions \cite{Mar}. This theorem occurred to be the milestone of powerful computer method for calculating integrals for many problems in differential and integral equations. The  package MATHEMATICA of  Wolfram Research is based on this theorem in calculating integrals.

 \theor\label{tnorm}
Buschman--Erdelyi operator of zero order smoothness $B_{0+}^{\nu,1}$ defined by (\ref{73})
acts under the Mellin transform as convolution (\ref{con}) with multiplicator
\begin{equation}\label{712}
m(s)=\frac{\Gamma(-s/2+\frac{\nu}{2}+1)\Gamma(-s/2-\frac{\nu}{2}+1/2)}
{\Gamma(1/2-\frac{s}{2})\Gamma(1-\frac{s}{2})}\label{s1}
\end{equation}
 for $\Re s<\min(2+\Re\nu,1-\Re\nu)$. Its norm is a periodic in $\nu$ and equals
\begin{equation}\label{s2}
\|B_{0+}^{\nu,1}\|_{L_2}=\frac{1}{\min(1,\sqrt{1-\sin\pi\nu})}.
\end{equation}
This operator is bounded in $L_2(0,\infty)$ if $\nu\neq 2k+1/2, k\in \mathbb{Z}$ and unbounded if
$\nu= 2k+1/2, k\in \mathbb{Z}$.

\vspace{3mm}

This theorem is the most important result of this article so we give a complete proof.

1. First let us prove the formula (\ref{712}) with a proper multiplicator. For it using consequently formulas   (7), p.~130, (2) p.~129, (4) p.~130 from \cite{Mar} we evaluate
$$
M(B_{0+}^{\nu,1})(s)=\frac{\Gamma(2-s)}{\Gamma(1-s)}M
\left[\int_0^\infty
\left\{H(\frac{x}{y}-1)P_\nu (\frac{x}{y})
\right\}
\left\{y f(y)\right\}\frac{dy}{y}
\right](s-1)=
$$
$$
=\frac{\Gamma(2-s)}{\Gamma(1-s)} M
\left[(x^2-1)_+^0P_\nu^0 (x)
\right]
(s-1)M\left[f\right](s),
$$
we use notations from \cite{Mar} for Heaviside and cutting power functions
$$
x_+^\alpha=\left\{
\begin{array}{rl}
x^\alpha, & \mbox{если } x\geqslant 0 \\
0, & \mbox{если } x<0 \\
\end{array}\right.
,\ H(x)=x_+^0=\left\{
\begin{array}{rl}
1, & \mbox{если } x\geqslant 0 \\
0, & \mbox{если } x<0. \\
\end{array}\right.
$$
Further using formulas 14(1) p.~234 и 4 p.~130 from \cite{Mar} we evaluate
$$
M\left[(x-1)_+^0 P_\nu^0 (\sqrt x)
\right](s)=
\frac{\Gamma(\frac{1}{2}+\frac{\nu}{2}-s)\Gamma(-\frac{\nu}{2}-s)}
{\Gamma(1-s)\Gamma(\frac{1}{2}-s)},
$$
$$
M\left[(x^2-1)_+^0 P_\nu^0 (x)
\right](s-1)=\frac{1}{2}\cdot
\frac
{
\Gamma(\frac{1}{2}+\frac{\nu}{2}-\frac{s-1}{2})
\Gamma(-\frac{\nu}{2}-\frac{s-1}{2})
}
{
\Gamma(1-\frac{s-1}{2})\Gamma(\frac{1}{2}-\frac{s-1}{2})
}=
$$
$$
=\frac{1}{2}\cdot\frac
{
\Gamma(-\frac{s}{2}+\frac{\nu}{2}+1)
\Gamma(-\frac{s}{2}-\frac{\nu}{2}+\frac{1}{2})
}
{\Gamma(-\frac{s}{2}+\frac{3}{2})\Gamma(-\frac{s}{2}+1)}
$$
under conditions $\Re s<\min(2+\Re\nu,1-\Re\nu)$. Now evaluate formula for
$$
M(B_{0+}^{\nu,1})(s)=\frac{1}{2}\cdot\frac{\Gamma(2-s)}{\Gamma(1-s)}\cdot
{\Gamma(-\frac{s}{2}+\frac{3}{2})\Gamma(-\frac{s}{2}+1)}.
$$
Applying to $\Gamma(2-s)$ the Legendre duplication formula (cf. \cite{BE1}) we evaluate
$$
M(B_{0+}^{\nu,1})(s)=\frac{2^{-s}}{\sqrt\pi}\cdot
\frac
{
\Gamma(-\frac{s}{2}+\frac{\nu}{2}+1)
\Gamma(-\frac{s}{2}-\frac{\nu}{2}+\frac{1}{2})
}
{\Gamma(1-s)}.
$$
Apply the Legendre duplication formula once more to $\Gamma(1-s)$ and the formula for  the multiplicator (\ref{s1}) is proved.
In the paper \cite{Sit3} it was shown that restrictions may be reduced to $0<\Re s<1$ for proper $\nu$. These restrictions may be weakened because they were derived for the class of all hypergeometric functions but we need just one special case of the Legendre function for which specified restrictions may be easily verified directly.

2. Now prove the formula (\ref{s2}) for a norm. From the multiplicator value we just found and theorem 4.7 from \cite{Sit1} on the line $\Re s=1/2, s=i u+1/2$ it follows
$$
|M(B_{0+}^{\nu,1})(i u+1/2)|=\frac{1}{\sqrt{2\pi}}\left|\frac
{
\Gamma(-i\frac{u}{2}-\frac{\nu}{2}+\frac{1}{4})
\Gamma(-i\frac{u}{2}+\frac{\nu}{2}+\frac{3}{4})
}
{\Gamma(\frac{1}{2}-iu)}\right|.
$$
Below operator symbol in  multiplicator will be omitted. Use  formulas for modulus $|z|=\sqrt{z\bar{z}}$
and gamma--function $\overline{\Gamma(z)}=\Gamma(\bar z)$ following from its definition as integral. The last property is true in general for the class of real--analytic functions. So we derive
$$
|M(B_{0+}^{\nu,1})(i u+1/2)|=
$$
$$
=\frac{1}{\sqrt{2\pi}}\left|\frac
{
\Gamma(-i\frac{u}{2}-\frac{\nu}{2}+\frac{1}{4})
\Gamma(i\frac{u}{2}-\frac{\nu}{2}+\frac{1}{4})
\Gamma(-i\frac{u}{2}+\frac{\nu}{2}+\frac{3}{4})
\Gamma(i\frac{u}{2}+\frac{\nu}{2}+\frac{3}{4})
}
{\Gamma(\frac{1}{2}-iu)\Gamma(\frac{1}{2}+iu)}\right|.
$$
In the numerator combine outer and inner terms and transform three pair of gamma-functions by the formula (см. \cite{BE1})
$$
\Gamma(\frac{1}{2}+z)\  \Gamma(\frac{1}{2}-z)=\frac{\pi}{\cos \pi z}.
$$
As a result we evaluate
$$
|M(B_{0+}^{\nu,1})(i u+1/2)|=
\sqrt{
\frac{\cos(\pi i u)}
{2\cos\pi(\frac{\nu}{2}+\frac{1}{4}+i\frac{u}{2})
\cos\pi(\frac{\nu}{2}+\frac{1}{4}-i\frac{u}{2})}
}=
$$
$$
=\sqrt{
\frac{\ch(\pi i u)}{\ch\pi u-\sin\pi\nu}
}
$$
Further denote as $t=\ch\pi u, 1\leqslant t <\infty$. So derive once more applying theorem 4.7 from \cite{Sit1}
$$
\sup_{u\in\mathbb{R}} |m(i u+\frac{1}{2})|=\sup_{1\leqslant t <\infty}
\sqrt{
\frac{t}{t-\sin\pi\nu}
}.
$$
So if $\sin\pi\nu\geqslant 0$ then supremum achieved at $t=1$ and for the norm the  formula (\ref{s2}) is valid
$$
\|B_{0+}^{\nu,1}\|_{L_2}=\frac{1}{\sqrt{1-\sin\pi\nu}}.
$$
Otherwise if  $\sin\pi\nu\leqslant 0$ then supremum achieved at $t\to\infty$ and the next formula is valid $$
\|B_{0+}^{\nu,1}\|_{L_2}=1.
$$
This part of the theorem is proved.

3. From the explicit values for norms and above cited theorem 4.7 from \cite{Sit1} follow conditions of boundedness or unboundedness and periodicity.
The theorem is completely proved.

Now proceed to finding multiplicators for all Buschman--Erdelyi operator of zero order smoothness.

\theor\label{tmult}
Buschman--Erdelyi operator of zero order smoothness acts under the Mellin transform as convolutions (\ref{con}). For their multiplicators the next formulas are valid

\begin{eqnarray}
& & m_{{_1S_{0+}^{\nu}}}(s)=\frac{\Gamma(-\frac{s}{2}+\frac{\nu}{2}+1) \Gamma(-\frac{s}{2}-\frac{\nu}{2}+\frac{1}{2})}{\Gamma(\frac{1}{2}-\frac{s}{2})\Gamma(1-\frac{s}{2})}= \label{3.11}; \\
& & =\frac{2^{-s}}{\sqrt{ \pi}} \frac{\Gamma(-\frac{s}{2}-\frac{\nu}{2}+\frac{1}{2}) \Gamma(-\frac{s}{2}+\frac{\nu}{2}+1)}{\Gamma(1-s)} , Re\, s < \min (2 + Re \, \nu, 1- Re\, \nu) \nonumber ; \\
& & m_{{_1P_{0+}^{\nu}}}(s)=\frac{\Gamma(\frac{1}{2}-\frac{s}{2})\Gamma(1-\frac{s}{2})}{\Gamma(-\frac{s}{2}+\frac{\nu}{2}+1) \Gamma(-\frac{s}{2}-\frac{\nu}{2}+\frac{1}{2})},~ Re\, s < 1; \label{3.12} \\
& & m_{{_1P_{-}^{\nu}}}(s)=\frac{\Gamma(\frac{s}{2}+\frac{\nu}{2}+1) \Gamma(\frac{s}{2}-\frac{\nu}{2})}{\Gamma(\frac{s}{2})\Gamma(\frac{s}{2}+\frac{1}{2})}, Re \, s > \max(Re \, \nu, -1-Re\, \nu); \label{3.13} \\
& & m_{{_1S_{-}^{\nu}}}(s)=\frac{\Gamma(\frac{s}{2})\Gamma(\frac{s}{2}+\frac{1}{2})}{\Gamma(\frac{s}{2}+\frac{\nu}{2}+\frac{1}{2}) \Gamma(\frac{s}{2}-\frac{\nu}{2})}, Re \, s >0 \label{3.14}
\end{eqnarray}

The next formulas are valid for norms of Buschman--Erdelyi operator of zero order smoothness in $L_2$:
\begin{eqnarray}
& & \| _1{S_{0+}^{\nu}} \| = \| _1{P_{-}^{\nu}}\|= 1/ \min(1, \sqrt{1- \sin \pi \nu}), \label{3.22} \\
& & \| _1{P_{0+}^{\nu}}\| = \| _1{S_{-}^{\nu}}\|= \max(1, \sqrt{1- \sin \pi \nu}). \label{3.23}
\end{eqnarray}

Similar results are proved in \cite{Sit2}--\cite{Sit3} for power weight spaces.

Corollary 1. Norms of operators \eqref{73} -- \eqref{733} are periodic in $\nu$ with period 2  $\|X^{\nu}\|=\|X^{\nu+2}\|$,  $X^{\nu}$ is any of operators \eqref{73} -- \eqref{733}.

Corollary 2. Norms of operators ${_1 S_{0+}^{\nu}}$, ${_1 P_{-}^{\nu}}$ are not bounded in total, every norm is greater or equals to 1. Norms are equal to 1 if $\sin \pi \nu \leq 0$.  Operators ${_1 S_{0+}^{\nu}}$, ${_1 P_{-}^{\nu}}$ are unbounded in $L_2$ if and only if $\sin \pi \nu = 1$ (or $\nu=(2k) + 1/2,~k \in \mathbb{Z}$).

Corollary 3. Norms of operators ${_1 P_{0+}^{\nu}}$, ${_1 S_{-}^{\nu}}$ are all bounded in $\nu$, every norm is not greater then $\sqrt{2}$. Norms are equal to 1 if $\sin \pi \nu \geq 0$.
Operators ${_1 P_{0+}^{\nu}}$, ${_1 S_{-}^{\nu}}$ are bounded in $L_2$ for all $\nu$. Maximum of norm equals $\sqrt 2 $ is achieved if and only if $\sin \pi \nu = -1$ (или $\nu= -1/2+(2k) ,~k \in \mathbb{Z}$).

The most important property of Buschman--Erdelyi operators of zero order smoothness is unitarity for integer $\nu$. It is just the case if interpret for these parameters the operator  $L_{\nu}$ as angular momentum operator in quantum mechanics.

\theor The  operators \eqref{73} -- \eqref{733} are unitary in  $L_2$ if and only if the parameter $\nu$ is an integer. In this case pairs of operators
$({_1 S_{0+}^{\nu}}$, ${_1 P_{-}^{\nu}})$ and  $({_1 S_{-}^{\nu}}$, ${_1 P_{0+}^{\nu}})$
are mutually inverse.

To formulate an interesting special case let us propose that operators \eqref{73} -- \eqref{733} act on functions permitting outer or inner differentiation in integrals, it is enough to suppose that $x f(x) \to 0$ for $x \to 0$. Then for  $\nu=1$
\begin{equation}\label {3.25}
{_1{P_{0+}^{1}}f=(I-H_1)f,~_1{S_{-}^{1}}f=(I-H_2)f,}
\end{equation}
and $H_1,~ H_2$ are famous Hardy operators,
\begin{equation}\label {3.26}
{H_1 f = \frac{1}{x} \int\limits_0^x f(y) dy,~H_2 f = \int\limits_x^{\infty}  \frac{f(y)}{y} dy,}
\end{equation}
$I$ is the identic operator.

Corollary 4. Operators \eqref{3.25} are unitary in  $L_2$ and mutually inverse. They are transmutations for a pair of differential operators  $d^2 / d x^2$ и $d^2 / d x^2 - 2/x^2$.

Unitarity of shifted Hardy operators \eqref{3.25} in $L_2$ is a known fact \cite{KuPe}. Below in application section we introduce a new class of generalizations for classical Hardy operators.

Now we list some properties of operators acting as convolutions by the formula \eqref{con} and with some multiplicator under the Mellin transform and being transmutations for the second derivative and angular momentum operator in quantum mechanics.

\theor\label{tOPmult} Let an operator $S_{\nu}$ acts by formulas \eqref{con}
and \eqref{76}. Then

а) its multiplicator satisfy a functional equation

\begin{equation}\label{5.1}{m(s)=m(s-2)\frac{(s-1)(s-2)}{(s-1)(s-2)-\nu(\nu+1)};}
\end{equation}

б) if any function $p(s)$ is periodic with period 2 ($p(s)=p(s-2)$) then a function $p(s)m(s)$ is a multiplicator for a new transmutation operator $S_2^{\nu}$ also acting by the rule \eqref{76}.

This theorem confirms the importance of studying transmutations in terms of the Mellin transform and
multiplicator functions.

Define the Stieltjes transform by (cf. \cite{KK})
$$
(S f)(x)= \int\limits_0^{\infty} \frac{f(t)}{x+t} dt.
$$
This operator also acts by the formula \eqref{con} with multiplicator $p(s)= \pi /sin (\pi s)$, it is bounded in $L_2$. Obviously  $p(s)=p(s-2)$. So from the theorem \ref{tOPmult} it follows that a convolution of the Stieltjes transform with bounded transmutations \eqref{73}--\eqref{733} are also transmutations of the same class bounded in $L_2$.

On this way many new classes of transmutations were introduced with special function kernels.

\section{\\Buschman--Erdelyi transmutations of the second kind.}

\setcounter{equation}{41}

Now we consider Buschman--Erdelyi transmutations of the second kind.

Definition 5.
Define a new pair of Buschman--Erdelyi transmutations of the second kind with Legendre functions of the second kind in kernels
\begin{equation}\label{6.1}
{{_2S^{\nu}}f=\frac{2}{\pi} \left( - \int\limits_0^x (x^2-y^2)^{-\frac{1}{2}}Q_{\nu}^1 (\frac{x}{y}) f(y) dy  +
\int\limits_x^{\infty} (y^2-x^2)^{-\frac{1}{2}}\mathbb{Q}_{\nu}^1 (\frac{x}{y}) f(y) dy\right),}
\end{equation}

\begin{equation}\label{6.2} {{_2P^{\nu}}f=\frac{2}{\pi} \left( - \int\limits_0^x (x^2-y^2)^{-\frac{1}{2}}\mathbb{Q}_{\nu}^1 (\frac{y}{x}) f(y) dy  -
\int\limits_x^{\infty} (y^2-x^2)^{-\frac{1}{2}}Q_{\nu}^1 (\frac{y}{x}) f(y) dy\right).}
\end{equation}

These operators are analogues of Buschman--Erdelyi transmutations of zero order smoothness.
If $y \to x \pm 0$ then integrals are defined by principal values. It is proved that they are transmutations of Sonine type for \eqref{6.1} and of Poisson type for \eqref{6.2}.

\theor  Operators  \eqref{6.1} -- \eqref{6.2} are of the form \eqref{con}
with multiplicators
\begin{eqnarray}
& & m_{_2S^{\nu}}(s)=p(s) \ m_{_1S_{-}^{\nu}}(s), \label{6.3}\\
& & m_{_2P^{\nu}}(s)=\frac{1}{p(s)} \ m_{_1P_{-}^{\nu}}(s), \label{6.4}
\end{eqnarray}
with  multiplicators of operators ${_1S_-^{\nu}}$, ${_1P_-^{\nu}}$ defined by \eqref{3.13} -- \eqref{3.14}, with period 2 function $p(s)$ equals
\begin{equation}\label {6.5}{p(s)=\frac{\sin \pi \nu+ \cos \pi s}{\sin \pi \nu - \sin \pi s}.}
\end{equation}

\theor The next formulas for norms are valid
\begin{eqnarray}
& &  \| {_2S^{\nu}} \|_{L_2}= \max (1, \sqrt {1+\sin \pi \nu}) , \label{6.9} \\
& &  \| {_2P^{\nu}} \|_{L_2}= 1 / {\min (1, \sqrt {1+\sin \pi \nu})} . \label{6.10}
\end{eqnarray}

Corollary.  Operator ${_2S^{\nu}}$ is bounded for all  $\nu$. Operator ${_2P^{\nu}}$ is not bounded if and only if then
$\sin \pi \nu=-1$.

\theor Operators ${_2S^{\nu}}$ and ${_2P^{\nu}}$ are unitary in $L_2$ if and only if $\nu\in\mathbb{Z}$.

\theor Let $\nu=i \beta+1/2,~\beta \in \mathbb{R}$. Then
\begin{equation}\label{6.11}{\| {_2S^{\nu}} \|_{L_2}=\sqrt{1+\ch \pi \beta},~\| {_2P^{\nu}} \|_{L_2}=1.}
\end{equation}

\theor The next formulas are valid
\begin{eqnarray}
& & {_2S^0} f = \frac{2}{\pi} \int\limits_0^{\infty} \frac{y}{x^2-y^2}f(y)\,dy, \label{6.12} \\
& & {_2S^{-1}} f = \frac{2}{\pi} \int\limits_0^{\infty} \frac{x}{x^2-y^2}f(y)\,dy. \label{6.13}
\end{eqnarray}

So in this case the operator ${_2S^{\nu}}$ reduce to a pair of semiaxis Hilbert transforms \cite{KK}.

For operators of the second kind also introduce more general ones with two parameters analogically to Buschman--Erdelyi transmutations of the first kind by formulas

\begin{eqnarray}
& & {_2S^{\nu,\mu}}f=\frac{2}{\pi} \left( \int\limits_0^x (x^2+y^2)^{-\frac{\mu}{2}} e^{-\mu \pi i} Q_{\nu}^{\mu}( \frac{x}{y}) f(y)\, dy\right. + \label{6.6} \\
& & +\left. \int\limits_x^{\infty} (y^2+x^2)^{-\frac{\mu}{2}}\mathbb{Q}_{\nu}^{\mu} (\frac{x}{y}) f(y)\, dy\right) \nonumber,
\end{eqnarray}
here $Q_{\nu}^{\mu}(z)$ is the Legendre function of the second kind, $\mathbb{Q}_{\nu}^{\mu}(z)$ is this function on the cut \cite{BE1}, $Re\, \nu < 1$. The second operator may be defined as formally conjugate in $L_2(0,\infty)$ to (\ref{6.6}).

\theor
The operator \eqref{6.6} on $C_0^{\infty}(0, \infty)$ is well defined and acts by
$$
M{[_2S^{\nu}]}(s)=m(s)\cdot M[x^{1-\mu} f](s), \label{6.7}
$$
\begin{eqnarray}
 m(s)=2^{\mu-1} \left( \frac{ \cos \pi(\mu-s) - \cos \pi \nu}{ \sin \pi(\mu-s) - \sin \pi \nu}  \right) \cdot\\
\cdot \left( \frac{\Gamma(\frac{s}{2})\Gamma(\frac{s}{2}+\frac{1}{2}))}{\Gamma(\frac{s}{2}+\frac{1-\nu-\mu}{2}) \Gamma(\frac{s}{2}+1+\frac{\nu-\mu}{2})} \right). \nonumber
\end{eqnarray}

\section{\\Buschman--Erdelyi transmutations of the third kind.}

\setcounter{equation}{53}

\subsection{Sonine--Katrakhov and Poisson--Katrakhov transmutations.}

Now we construct transmutations which are unitary for all $\nu$. They are defined by formulas
\begin{eqnarray}
& & S_U^{\nu} f = - \sin \frac{\pi \nu}{2}\  {_2S^{\nu}}f+ \cos \frac{\pi \nu}{2}\  {_1S_-^{\nu}}f, \label{6.14} \\
& & P_U^{\nu} f = - \sin \frac{\pi \nu}{2}\  {_2P^{\nu}}f+ \cos \frac{\pi \nu}{2}\  {_1P_-^{\nu}}f. \label{6.15}
\end{eqnarray}
For all values $\nu \in \mathbb{R}$ they are linear combinations of Buschman--Erdelyi transmutations of the first and second kinds of zero order smoothness. Also they are in the defined below class of Buschman--Erdelyi transmutations of the third kind.
Integral representations are valid

\begin{eqnarray}
& & S_U^{\nu} f = \cos \frac{\pi \nu}{2} \left(- \frac{d}{dx} \right) \int\limits_x^{\infty} P_{\nu}\lr{\frac{x}{y}} f(y)\,dy +  \label{6.16}\\
& & + \frac{2}{\pi} \sin \frac{\pi \nu}{2} \left(  \int\limits_0^x (x^2-y^2)^{-\frac{1}{2}}Q_{\nu}^1 \lr{\frac{x}{y}} f(y)\,dy  \right.-
%\nonumber\\& & -
 \int\limits_x^{\infty} (y^2-x^2)^{-\frac{1}{2}}\mathbb{Q}_{\nu}^1 \lr{\frac{x}{y}} f(y)\,dy \Biggl. \Biggr), \nonumber \\
& & P_U^{\nu} f = \cos \frac{\pi \nu}{2}  \int\limits_0^{x} P_{\nu}\lr{\frac{y}{x}} \left( \frac{d}{dy} \right) f(y)\,dy - \label{6.17} \\
& &  -\frac{2}{\pi} \sin \frac{\pi \nu}{2} \left( - \int\limits_0^x (x^2-y^2)^{-\frac{1}{2}}\mathbb{Q}_{\nu}^1\lr{\frac{y}{x}} f(y)\,dy   \right.-
%\nonumber\\& & -
\int\limits_x^{\infty} (y^2-x^2)^{-\frac{1}{2}} Q_{\nu}^1 \lr{\frac{y}{x}} f(y)\,dy \Biggl. \Biggr). \nonumber
\end{eqnarray}

\theor Operators \eqref{6.14}--\eqref{6.15},  \eqref{6.16}--\eqref{6.17} for all $\nu\in\mathbb{R}$
are unitary, mutually inverse and conjugate in $L_2$. They are transmutations acting by \eqref{75}.
$S_U^{\nu}$ is a Sonine type transmutation and $P_U^{\nu}$ is a Poisson type one.

Transmutations like (\ref{6.16})--(\ref{6.17}) but with kernels into more complicated form with hypergeometric functions were first introduced by Katrakhov in 1980. Due to it the author propose terms for this class of operators as Sonine--Katrakhov and Poisson--Katrakhov. In author's papers these operators were reduced to more simple form of Buschman--Erdelyi ones. It made possible to include this class of operators in general composition (or factorization) method  \cite{SiKa2}, \cite{SiKa3}, \cite{Sit5}.

\subsection{Buschman--Erdelyi transmutations
of the third kind with arbitrary weight function.}

Define sine and cosine Fourier transforms with inverses

\begin{equation}\label{1}{F_c f = \sqrt{\frac{2}{\pi}} \int\limits_0^{\infty} f(y) \cos (t y) \, dy,
\ \ F_c^{-1}=F_c,}
\end{equation}

\begin{equation}\label{2}{F_s f = \sqrt{\frac{2}{\pi}} \int\limits_0^{\infty} f(y) \sin (t y) \, dy,}
\ \ F_s^{-1}=F_s.
\end{equation}

Define Hankel (Fourier--Bessel) transform and its inverse by

\begin{eqnarray}
& &  F_{\nu} f = \frac{1}{2^{\nu} \Gamma (\nu+1 )} \int\limits_0^{\infty} f(y)\, j_{\nu}(t y) \, y^{2 \nu + 1} \, dy = \nonumber\\
& & = \int\limits_0^{\infty} f(y) \frac{J_\nu(t y)}{(t y)^{\nu}}\,  y^{2 \nu + 1} \, dy = \frac{1}{t^{\nu}} \int\limits_0^{\infty} f(y)  J_{\nu}(t y) \, y^{\nu + 1} \, dy, \label{3} \\
& & F_{\nu}^{-1} f = \frac{1}{(y)^{\nu}} \int\limits_0^{\infty} f(t) J_{\nu}(y t)\, t^{\nu + 1} \, dt. \label{4}
\end{eqnarray}
Here $J_\nu(\cdot)$ is the Bessel function \cite{BE1}, \  $j_\nu(\cdot)$ is normalized Bessel function \cite{Kip1}.
Operators \eqref{1}-\eqref{2} are unitary self--conjugate in $L_2(0, \infty)$.  Operators \eqref{3}--\eqref{4} are unitary self--conjugate in power weighted space $L_{2,\, \nu}(0, \infty)$.

Now define on proper functions the first pair of Buschman--Erdelyi transmutations
of the third kind

\begin{eqnarray}
& &  S_{\nu,\, c}^{(\varphi)} = F^{-1}_c \left( \frac{1}{\varphi (t)} F_{\nu} \right), \label{5} \\
& &  P_{\nu,\, c}^{(\varphi)} = F^{-1}_{\nu} \left( \varphi (t) F_{c} \right), \label{6}
\end{eqnarray}

and the second pair by

\begin{eqnarray}
& &  S_{\nu,\, s}^{(\varphi)} = F^{-1}_s \left( \frac{1}{\varphi (t)} F_{\nu} \right), \label{7} \\
& &  P_{\nu,\, s}^{(\varphi)} = F^{-1}_{\nu} \left( \varphi (t) F_{s} \right). \label{8},
\end{eqnarray}
with $\varphi (t)$ being an arbitrary weight function.

The operators defined on proper functions are transmutations for $B_{\nu}$ and $D^2$. They may be expressed in the integral form.

\theor Define transmutations for $B_{\nu}$ and $D^2$ by formulas

$$
S^{(\varphi)}_{\nu, \left\{ \begin{matrix} s \\ c \end{matrix} \right\} }=
F^{-1}_{ \left\{ \begin{matrix} s \\ c \end{matrix} \right\} } \left( \frac{1}{\varphi(t)}F_{\nu} \right),
$$

$$
P^{(\varphi)}_{\nu, \left\{ \begin{matrix} s \\ c \end{matrix} \right\} }=
F^{-1}_{ \nu } \left( \varphi(t)F_{ \left\{ \begin{matrix} s \\ c \end{matrix} \right\} } \right).
$$

Then for the Sonine type transmutation an integral form is valid

\begin{equation}{\left( S^{(\varphi)}_{\nu, \left\{ \begin{matrix} s \\ c \end{matrix} \right\} } f \right) (x) = \sqrt{\frac{2}{\pi}} \int\limits_0^{\infty} K(x, y) f(y) \, dy,}
\end{equation}
где
$$
K(x, y) = y^{\nu + 1} \int\limits_0^{\infty} \frac{\left\{ \begin{matrix} \sin (x t) \\ \cos (x t) \end{matrix} \right\}}{\varphi (t)\, t^{\nu} } J_{\nu} (y t) dt.
$$

For the Poisson type transmutation an integral form is valid

\begin{equation}{\left( P^{(\varphi)}_{\nu, \left\{ \begin{matrix} s \\ c \end{matrix} \right\} } f \right) (x) = \sqrt{\frac{2}{\pi}} \int\limits_0^{\infty} G(x, y) f(y) \, dy,}
\end{equation}
где
$$
G(x, y) = \frac{1}{x^\nu} \int\limits_0^{\infty} \varphi (t) \, t^{\nu + 1} \left\{ \begin{matrix} \sin (y t) \\ \cos (y t) \end{matrix}\right\}  J_{\nu} (x t) dt.
$$

Introduced before unitary transmutations of Sonine--Katrakhov and Poisson--Katrakhov are special cases of this class operators. For this case we must choice a weight function $\varphi (t)$ as a power function depending on the parameter $\nu$. The author plans to publish a special paper on Buschman--Erdelyi transmutations
of the third kind with arbitrary weight function.

\section{\\Some Applications of Buschman--Erdelyi transmutations.}

\setcounter{equation}{57}

In this section we gather some applications of Buschman--Erdelyi operators (but not all). Due to the article restrictions most items are only  briefly mentioned with most informative facts and instructive references. Some applications only mentioned as problems for future investigations.

\bigskip

\subsection{Norm estimates and embedding theorems in Kipriyanov spaces.}
\medskip

Consider a set of functions $\mathbb{D}(0, \infty)$. If $f(x) \in \mathbb{D}(0, \infty)$ then $f(x) \in C^{\infty}(0, \infty),~ f(x)$ is zero at infinity. On this set define seminorms
\begin{eqnarray}
& & \|f\|_{h_2^{\alpha}}=\|D_{-}^{\alpha}f\|_{L_2(0, \infty)} \label{7.1} \\
& &  \|f\|_{\widehat{h}_2^{\alpha}}=\|x^{\alpha} (-\frac{1}{x}\frac{d}{dx})^{\alpha}f\|_{L_2(0, \infty)} \label{7.2}
\end{eqnarray}
here $D_-^{\alpha}$ is the Riemann--Liouville fractional integrodifferentiation, operator in  \eqref{7.2} is defined by
\begin{equation}\label{7.3}{(-\frac{1}{x}\frac{d}{dx})^{\beta}=2^{\beta}I_{-; \, 2, \,0}^{-\beta}x^{-2 \beta},}
\end{equation}
$I_{-; 2, \, 0}^{-\beta}$ is Erdelyi--Kober operator, $\alpha\in\mathbb{R}$.  For $\beta = n \in \mathbb{N}_0$ expression  \eqref{7.3} reduces to classical derivatives.

\theor
Let $f(x) \in \mathbb{D}(0, \infty)$. Then the next formulas are valid:
\begin{eqnarray}
& & D_{-}^{\alpha}f={_1S_{-}^{\alpha-1}} {x^{\alpha} (-\frac{1}{x}\frac{d}{dx})^{\alpha}} f, \label{7.4} \\
& & x^{\alpha} (-\frac{1}{x}\frac{d}{dx})^{\alpha}f={_1P_{-}^{\alpha-1}} D_{-}^{\alpha}f. \label{7.5}
\end{eqnarray}

So Buschman--Erdelyi transmutations of zero order smoothness for $\alpha \in \mathbb{N}$ links differential operators in seminorms definitions \eqref{7.1} and \eqref{7.2}.

\theor
Let $f(x) \in \mathbb{D}(0, \infty)$. Then the next inequalities hold true for seminorms
\begin{eqnarray}
& &  \|f\|_{h_2^{\alpha}} \leq \max (1, \sqrt{1+\sin \pi \alpha}) \|f\|_{\widehat{h}_2^{\alpha}}, \label{7.7}\\
& & \|f\|_{\widehat{h}_2^{\alpha}} \leq \frac{1}{\min (1, \sqrt{1+\sin \pi \alpha})} \|f\|_{h_2^{\alpha}}, \label{7.8}
\end{eqnarray}
here $\alpha$ is any real number except $\alpha \neq -\frac{1}{2}+2k,~k \in \mathbb{Z}$.

Constants in inequalities \eqref{7.7}--\eqref{7.8} are not greater than 1, it will be used below.
If $\sin \pi \alpha = -1 $ or  $\alpha = -\frac{1}{2}+2k,~k \in \mathbb{Z}$ then the estimate \eqref{7.8} is not valid.

Define on  $\mathbb{D} (0, \infty )$ the Sobolev norm
\begin{equation}\label{7.9}{\|f\|_{W_2^{\alpha}}=\|f\|_{L_2 (0, \infty)}+\|f\|_{h_2^{\alpha}}.}
\end{equation}
Define one more norm
\begin{equation}\label{7.10}{\|f\|_{\widehat{W}_2^{\alpha}}=\|f\|_{L_2 (0, \infty)}+\|f\|_{\widehat{h}_2^{\alpha}}}
\end{equation}
Define spaces $W_2^{\alpha},~ \widehat{W}_2^{\alpha}$ as closures of $D(0,
\infty)$ in  \eqref{7.9} or \eqref{7.10} respectively.

\theor а) For all $\alpha \in \mathbb{R}$ the space $\widehat{W}_2^{\alpha}$ is continuously imbedded in  $W_2^{\alpha}$, moreover
\begin{equation}\label{7.11}{\|f\|_{W_2^{\alpha}}\leq A_1 \|f\|_{\widehat{W}_2^{\alpha}},}
\end{equation}
with $A_1=\max (1, \sqrt{1+\sin \pi \alpha})$.

б) Let $\sin \pi \alpha \neq -1$ or $\alpha \neq -\frac{1}{2} + 2k, ~ k \in \mathbb{Z}.  $ Then the inverse embedding of $W_2^{\alpha}$  in $\widehat{W}_2^{\alpha}$ is valid, moreover
\begin{equation}\label{7.12}{\|f\|_{\widehat{W}_2^{\alpha}}\leq A_2 \|f\|_{W_2^{\alpha}},}
\end{equation}
with $A_2 =1/  \min (1, \sqrt{1+\sin \pi \alpha})$.

в) Let $\sin \pi \alpha \neq -1$, then spaces $W_2^{\alpha}$  and $\widehat{W}_2^{\alpha}$ are isomorphic with equivalent norms.

г) Constants in embedding inequalities \eqref{7.11}--\eqref{7.12} are  sharp.

In fact this theorem is a direct corollary of results on boundedness and norm estimates in $L_2$ of
Buschman--Erdelyi transmutations of zero order smoothness. At the same manner from unitarity of these operators follows the next

\theor Norms
\begin{eqnarray}
& & \|f\|_{W_2^{\alpha}} = \sum\limits_{j=0}^s \| D_{-}^j f\|_{L_2}, \label{7.13} \\
& & \|f\|_{\widehat{W}_2^{\alpha}}=\sum\limits_{j=0}^s \| x^j(-\frac{1}{x}\frac{d}{dx})^j f \|_{L_2} \label{7.14}
\end{eqnarray}
are equivalent   for integer  $s \in \mathbb{Z}$. Moreover every term in  \eqref{7.13} equals to appropriate term in \eqref{7.14} of the same index $j$.

I.\,Kipriyanov introduced in \cite{Kip2} function spaces which essentially influenced the theory of partial differential equations with Bessel operators and in more general sense theory of singular and degenerate equations. These spaces are defined by the next way. First consider subset of even functions in  $\mathbb{D}(0, \infty)$ with all zero derivatives of odd orders at  $x=0$. Denote this set as $\mathbb{D}_c (0, \infty)$ and equipped it with a norm
\begin{equation}\label{7.15}{\|f\|_{\widetilde{W}_{2, k}^s} = \|f\|_{L_{2, k}}+\|B_k^{\frac{s}{2}}\|_{L_{2, k}}}
\end{equation}
here $s$ is an even natural number,  $B^{s/2}_k$ is an iteration of the Bessel operator. Define Kipriyanov spaces for even  $s$ as a closure of $D_c (0, \infty)$ in the norm \eqref{7.15}. It is a known fact that equivalent to \eqref{7.15} norm may be defined by \cite{Kip2}
\begin{equation}\label{7.16}{\|f\|_{\widetilde{W}_{2, k}^s} = \|f\|_{L_{2, k}}+\|x^s(-\frac{1}{x}\frac{d}{dx})^s f\|_{L_{2, k}}}
\end{equation}
So the norm $\widetilde{W}_{2, \, k}^s$ may be defined for all $s$. Essentially this approach is the same as in   \cite{Kip2}, another approach is based on usage of Hankel transform. Below we adopt the norm \eqref{7.16} for the space $\widetilde{W}_{2, k}^s$.

Define weighted Sobolev norm by
\begin{equation}\label{7.17}{\|f\|_{W_{2, k}^s} = \|f\|_{L_{2, k}}+\|D_{-}^s f\|_{L_{2, k}}}
\end{equation}
and a space $W_{2, \, k}^s$ as a closure of $\mathbb{D}_c (0, \infty)$ in this norm.

\theor\label{tvloz1}
а) Let $k \neq -n, ~ n \in \mathbb{N}$. Then the space  $\widetilde{W}_{2, \, k}^s$ is continuously embedded into  $W_{2, \, k}^s$, and there exist a constant $A_3>0$ such that
\begin{equation}\label{7.18}{\|f\|_{W_{2, k}^s}\leq A_3 \|f\|_{\widetilde{W}_{2, k}^s},}
\end{equation}
б) Let $k+s \neq -2m_1-1, ~ k-s \neq -2m_2-2, ~ m_1 \in \mathbb{N}_0, ~ m_2 \in \mathbb{N}_0$. Then the inverse embedding holds true of $W_{2, \, k}^s$ into $\widetilde{W}_{2, \, k}^s$, and there exist a constant $A_4>0$, such that
\begin{equation}\label{7.19}{\|f\|_{\widetilde{W}_{2, k}^s}\leq A_4 \|f\|_{W_{2, k}^s}.}
\end{equation}
в) If the above mentioned conditions are not valid then embedding theorems under considerations fail.

Corollary 1.  Let the next conditions hold true: $k \neq -n, ~ n \in \mathbb{N}$; $k+s \neq -2m_1-1,  ~ m_1 \in \mathbb{N}_0; ~ k-s \neq -2m_2-2, ~ m_2 \in \mathbb{N}_0$. Then Kipriyanov spaces may be defined as closure of $D_c (0, \infty)$ in the weighted Sobolev norm \eqref{7.17}.

Corollary 2.  Sharp constants in embedding theorems \eqref{7.18}--\eqref{7.19} are
$$
A_3 = \max (1, \|{_1S_-^{s-1}} \| _ {L_{2, k}}), ~ A_4=\max(1, \|{_1P_-^{s-1}}\|_{L_{2, k}}).
$$

It is obvious that the theorem above and its corollaries are direct consequences of estimates for
Buschman--Erdelyi transmutations. Sharp constants in embedding theorems \eqref{7.18}--\eqref{7.19} are also direct  consequences of estimates for
Buschman--Erdelyi transmutations of zero order smoothness. Estimates in  $L_{p, \alpha}$ not included in this article allow to consider embedding theorems for general Sobolev and Kipriyanov spaces.

So by applying Buschman--Erdelyi transmutations of zero order smoothness we received an answer to a problem which for a long time was discussed in "folklore": --- are Kipriyanov spaces isomorphic to power weighted Sobolev spaces or not? Of course we investigated just the simplest case, results may be generalize to other seminorms, higher dimensions, bounded domains but the principal idea is clear. All that do not in any sense  disparage neither essential role nor necessity for applications of Kipriyanov spaces in the theory of partial differential equations.

The importance of Kipriyanov spaces is a special case of the next general principle of L.\,Kudryavtsev:

\begin{center}
 "\textit{EVERY EQUATION MUST BE INVESTIGATED IN ITS OWN SPACE!}"
\end{center}

The proved in this section embedding theorems may be applied to direct transfer of known solution estimates for  $B$--elliptic equations in Kipriyanov spaces (cf. \cite{Kip1},\cite{Kip2} ) to new estimates in weighted Sobolev spaces, it is a direct consequence of boundedness and transmutation properties of
Buschman--Erdelyi transmutations.

\bigskip

\subsection{Solution representations to partial differential equations with Bessel operators.}

The above classes of transmutations may be used for deriving explicit formulas for solutions of partial differential equations with Bessel operators via unperturbed equation solutions. An example is the  $B$--elliptic equation of the form
\begin{equation}
\sum_{k=1}^{n}B_{\nu,x_k}u(x_1,\dots, x_n)=f,
\end{equation}
and similar  $B$--hyperbolic and $B$--parabolic equations. This idea early works by Sonine--Poisson--Delsarte transmutations, cf. \cite{Car1}--\cite{Car3}, \cite{CarSho}, \cite{Kip1}.
New results follow automatically for new classes of transmutations.

\bigskip

\subsection{Cauchy problem for Euler--Poisson--Delsarte equation (EPD).}

\medskip

Consider EPD equation in a half space
$$
B_{\alpha,\, t} u(t,x)= \frac{ \partial^2 \, u }{\partial t^2} + \frac{2 \alpha+1}{t} \frac{\partial u}{\partial t}=\Delta_x u+F(t, x),
$$
with $t>0,~x \in \mathbb{R}^n$.
Let us consider a general plan for finding different initial and boundary conditions at  $t=0$ with guaranteed existence of solutions. Define any transmutations $X_{\alpha, \, t}$ and $Y_{\alpha, \, t}$
satisfying  \eqref{75}. Suppose that functions $X_{\alpha, \, t} u=v(t,x)$, $X_{\alpha, \, t} F=G(t,x)$ exist. Suppose that unperturbed Cauchy problem
\begin{equation}\label{7.28}{\frac{ \partial^2 \, v }{\partial t^2} =\Delta_x v+G,~ v|_{t=0}=\varphi (x),~ v'_t|_{t=0}=\psi (x)}
\end{equation}
is correctly solvable in a half space. Then if $Y_{\alpha, \, t}=X^{-1}_{\alpha, \, t}$
then we receive the next initial conditions
\begin{equation}\label{7.29}{X_{\alpha} u|_{t=0}=a(x),~(X_{\alpha} u)'|_{t=0}=b(x).}
\end{equation}
By this method   the choice of different classes of transmutations (Sonine--Poisson--Delsarte, Buschman--Erdelyi of the first, second and third kinds,  Buschman--Erdelyi of the zero order smoothness, unitary transmutations of Sonine--Katrakhov and Poisson--Katrakhov, transmutations with general kernels) will correspond different kinds of initial conditions \cite{Sit3}.

In the monograph of Pskhu \cite{Psh} this method is applied for solving an equation with fractional derivatives with the usage of Stankovic transform. Glushak applied Buschman--Erdelyi operators in \cite{Glu}.

The Buschman--Erdelyi operators were first introduced exactly for EPD equation by Copson. We formulate his result now.

\bigskip

Copson lemma.
\medskip

Consider partial differential equation with two variables on the plane

$$
\frac{\partial^2 u(x,y)}{\partial x^2}+\frac{2\alpha}{x}\frac{\partial u(x,y)}{\partial x}=
\frac{\partial^2 u(x,y)}{\partial y^2}+\frac{2\beta}{y}\frac{\partial u(x,y)}{\partial y}
$$
(this is EPD equation or $B$--hyperbolic one in Kipriyanov's terminology)
for $x>0,\  y>0$ and $\beta>\alpha>0$
with boundary conditions on characteristics

$$
u(x,0)=f(x), u(0,y)=g(y), f(0)=g(0).
$$

It is supposed that the solution $u(x,y)$ is continuously differentiable in the closed first quadrant and has second derivatives in this open  quadrant, boundary functions $f(x), g(y)$ are differentiable.

Then if the solution exist the next formulas hold true

\begin{equation}
\label{Cop1}
\frac{\partial u}{\partial y}=0, y=0,  \frac{\partial u}{\partial x}=0, x=0,
\end{equation}
\begin{equation}
\label{Cop2}
2^\beta \Gamma(\beta+\frac{1}{2})\int_0^1 f(xt)t^{\alpha+\beta+1}
\lr{1-t^2}^{\frac{\beta -1}{2}}P_{-\alpha}^{1-\beta}{t}\,dt=
\end{equation}
\begin{equation*}
=2^\alpha \Gamma(\alpha+\frac{1}{2})\int_0^1 g(xt)t^{\alpha+\beta+1}
\lr{1-t^2}^{\frac{\alpha -1}{2}}P_{-\beta}^{1-\alpha}{t}\,dt,
\end{equation*}
$$
\Downarrow
$$
\begin{equation}
\label{Cop3}
g(y)=\frac{2\Gamma(\beta+\frac{1}{2})}{\Gamma(\alpha+\frac{1}{2})
\Gamma(\beta-\alpha)}y^{1-2\beta}
\int_0^y x^{2\alpha-1}f(x)
\lr{y^2-x^2}^{\beta-\alpha-1}x \,dx,
\end{equation}
here  $P_\nu^\mu(z)$ is the Legendre function of the first kind \cite{Sit1}.

So the main conclusion from Copson lemma is that data on characteristics can not be taken arbitrary, these functions must be connected by Buschman--Erdelyi operators of the first kind, for more detailed consideration cf. \cite{Sit1}.

\bigskip

\subsection{Applications to generalized translations.}
\medskip

This class of operators was thoroughly studied by Levitan  \cite{Lev1}--\cite{Lev2}.
It has many applications to partial differential equations, including Bessel operators \cite{Lev3}. Generalized translations are used for moving singular point from the origin to any location. They are explicitly expressed via transmutations \cite{Lev3}. Due to this fact new classes of transmutations lead to new classed of generalized translations.

\bigskip

\subsection{Applications to  Dunkl operators.}
\medskip

In recent years  Dunkl operators were thoroughly studied. These are difference--differentiation operators consisting of combinations of classical derivatives and finite differences. In higher dimensions Dunkl operators are defined by symmetry and reflection groups. For this class there are many results on transmutations which are of Sonine--Poisson--Delsarte and Buschman--Erdelyi types, cf. \cite{Sig} and references therein.

\bigskip

\subsection{Applications of Buschman--Erdelyi operators to the Radon transform.}
\medskip

It was proved by Ludwig in \cite{Lud} that the Radon transform in terms of spherical harmonics acts in every harmonics at radial components as Buschman--Erdelyi operators. Let us formulate this result.

\theor Ludwig theorem \ (\cite{Lud},\cite{Hel}). Let the function $f(x)$ being expanded in $\mathbb{R}^n$ by spherical harmonics
\begin{equation}
f(x)=\sum_{k,l}f_{k,l}(r) Y_{k,l}(\theta).
\end{equation}
Then the Radon transform of this function may be calculated as another series in spherical harmonics
\begin{equation}
Rf(x)=g(r,\theta)=\sum_{k,l}g_{k,l}(r) Y_{k,l}(\theta),
\end{equation}
\begin{equation}
\label{lu1}
g_{k,l}(r)=с(n)\int_r^\infty \lr{1-\frac{s^2}{r^2}}^{\frac{n-3}{2}} C_l^{\frac{n-2}{2}}\lr{\frac{s}{r}}
f_{k,l}(r) r^{n-2}\,ds,
\end{equation}
here $с(n)$ is some known constant, $C_l^{\frac{n-2}{2}}\lr{\frac{s}{r}}$ is the Gegenbauer function \cite{BE1}. The inverse formula is also valid of representing values $f_{k,l}(r)$ via $g_{k,l}(r)$.

The Gegenbauer function may be easily reduced to the Legendre function  \cite{BE1}. So the Ludwig's formula (\ref{lu1}) reduce the Radon transform in terms of spherical harmonics series and up to unimportant power and constant terms to Buschman--Erdelyi operators of the first kind.

Exactly this formula in dimension two was developed by Cormack as the first step to the Nobel prize.
Special cases of Ludwig's formula proved in 1966 are  for any special spherical harmonics and in the simplest case on pure radial function, in this case it is reduced to Sonine--Poisson--Delsarte transmutations of Erdelyi--Kober type. Besides the fact that such formulas are known  for about half a century they are   rediscovered  still... As consequences of the above connections the results may be proved for integral representations, norm estimates, inversion formulas for the Radon transform via
Buschman--Erdelyi operators. In particular it makes clear that different kinds of inversion formulas for the Radon transform are at the same time inversion formulas for the Buschman--Erdelyi transmutations of the first kind and vice versa. A useful reference for this approach is \cite{Dea}.

\bigskip

\subsection{ Application of Buschman--Erdelyi operators to generalized polynomials and spherical harmonics.}

It was known for many years that a problem of describing polynomial solutions for $B$--elliptic equation do not need the new theory. The answer is in the transmutation theory. A simple fact that Sonine--Poisson--Delsarte transmutations transform power function into another power function means that they also transform explicitly so called $B$--harmonic polynomials into classical harmonic polynomials and vice versa. The same is true for generalized  $B$--harmonics because they are restrictions of $B$--harmonic polynomials onto the unit sphere. This approach is thoroughly applied by Rubin  \cite{Rub1}--\cite{Rub2}.
Usage of Buschman--Erdelyi operators refresh this theory with new possibilities.

\bigskip

\subsection{Application of Buschman--Erdelyi transmutations for estimation
of generalized Hardy operators.}

We proved unitarity of shifted Hardy operators (\ref{3.25}) and mentioned that it is a known fact from \cite{KuPe}. It is interesting that Hardy operators naturally arise in transmutation theory. Use the theorem 7 with integer parameter which guarantees unitarity for finding more unitary in $L_2(0,\infty)$ integral operators of very simple form.

\theor The next are  pair of unitary mutually inverse integral operators in $L_2(0,\infty)$:
\begin{eqnarray*}
\label{84}
U_3f= f+\int_0^x f(y)\,\frac{dy}{y},\  U_4f= f+\frac{1}{x}\int_x^\infty f(y)\,dy,\\\nonumber
U_5f= f+3x\int_0^x f(y)\,\frac{dy}{y^2},\  U_6f= f-\frac{3}{x^2}\int_0^x y f(y)\,dy,\\\nonumber
U_7f=f+\frac{3}{x^2}\int_x^\infty y f(y)\,dy,\  U_8f=f-3x \int_x^\infty f(y)\frac{dy}{y^2},\\\nonumber
U_9f=f+\frac{1}{2}\int_0^x \left(\frac{15x^2}{y^3}-\frac{3}{y}\right)f(y)\,dy,\\\nonumber
U_{10}f=f+\frac{1}{2}\int_x^\infty \left(\frac{15y^2}{x^3}-\frac{3}{x}\right)f(y)\,dy.\\\nonumber
\end{eqnarray*}

\bigskip

\subsection{Integral operators with more general functions as kernels.}
\medskip

Consider an operator ${{}_1S_{0+}^{\nu}}$. It has the form
\begin{equation}\label{7.30}{{{_1}S_{0+}^{\nu}}=\frac{d}{dx}\int\limits_0^x K\lr{\frac{x}{y}} f(y) \, dy,}
\end{equation}
with kernel $K$ expressed by $K(z)=P_{\nu}(z)$. Simple properties of special functions lead to the fact that ${{_1}S_{0+}^{\nu}}$ is a special case of \eqref{7.30} with Gegenbauer function kernel
\begin{equation}\label{7.31}{K(z)=\frac{\Gamma(\alpha+1)\ \Gamma(2 \beta)}{2^{p-\frac{1}{2}}\Gamma(\alpha+2\beta)\ \Gamma(\beta+ \frac{1}{2})}(z^{\alpha}-1)^{\beta-\frac{1}{2}}C^{\beta}_{\alpha}(z)}
\end{equation}
with $\alpha=\nu, ~ \beta = \frac{1}{2}$ or with Jacobi function kernel
\begin{equation}\label{7.32}{K(z)=\frac{\Gamma(\alpha+1)}{2^{\rho}\Gamma(\alpha+\rho+1)}(z-1)^{\rho}(z+1)^{\sigma}P^{(\rho, \sigma)}_{\alpha}(z)}
\end{equation}
with $\alpha=\nu, ~ \rho=\sigma = 0$. More general are operators with Gauss hypergeometric function kernel ${{_2}F_1}$, Mayer  $G$ or Fox $H$ function kernels, cf. \cite{KK}, \cite{KiSa}. For studying such operators inequalities for kernel functions are very useful, e.g.  \cite{SiKa1}--\cite{SiKa2}.

Define the first class of generalized operators.

Definition 6.{ Define Gauss--Buschman--Erdelyi operators by formulas
\begin{equation}\label{7.33}{{_1F_{0+}}(a, b, c)[f]=\frac{1}{2^{c-1}\Gamma(c)}.}
\end{equation}
$$
\int\limits_0^x\lr{\frac{x}{y}-1}^{c-1}\lr{\frac{x}{y}+1}^{a+b-c} {{_2}F_1}\lr{^{a,b}_c| \frac{1}{2}-\frac{1}{2}\frac{x}{y}} f(y) \, dy,
$$
\begin{equation}\label{7.34}{{{_2}F_{0+}}(a, b, c)[f]=\frac{1}{2^{c-1}\Gamma(c)}.}
\end{equation}
$$
\int\limits_0^x\lr{\frac{y}{x}-1}^{c-1}\lr{\frac{y}{x}+1}^{a+b-c} {{_2}F_1}\lr{^{a,b}_c| \frac{1}{2}-\frac{1}{2}\frac{y}{x}} f(y) \, dy,
$$
\begin{equation}\label{7.35}{{_1F_{-}}(a, b, c)[f]=\frac{1}{2^{c-1}\Gamma(c)}.}
\end{equation}
$$
\int\limits_0^x\lr{\frac{y}{x}-1}^{c-1}\lr{\frac{y}{x}+1}^{a+b-c} {{_2}F_1}\lr{^{a,b}_c| \frac{1}{2}-\frac{1}{2}\frac{y}{x}} f(y) \, dy,
$$
\begin{equation}\label{7.36}{{{_2}F_{-}}(a, b, c)[f]=\frac{1}{2^{c-1}\Gamma(c)}.}
\end{equation}
$$
\int\limits_0^x\lr{\frac{x}{y}-1}^{c-1}\lr{\frac{x}{y}+1}^{a+b-c} {{_2}F_1}\lr{^{a,b}_c| \frac{1}{2}-\frac{1}{2}\frac{x}{y}} f(y) \, dy,
$$
\begin{eqnarray}
& {{_3}F_{0+}}[f]=\frac{d}{dx} {{_1}F_{0+}}[f], & {{_4}F_{0+}}[f]= {{_2}F_{0+}} \frac{d}{dx} [f], \label{7.37} \\
& {{_3}F_{-}}[f]={{_1}F_{-}}(-\frac{d}{dx} ) [f], & {{_4}F_{-}}[f]= (- \frac{d}{dx} ) {{_2}F_{-}}[f]. \label{7.38}
\end{eqnarray}
}

Symbol ${{_2}F_1}$ in definitions \eqref{7.34} and \eqref{7.36} means Gauss hypergeometric function on natural domain and in \eqref{7.33} and \eqref{7.35} the main branch of its analytical continuation.

Operators \eqref{7.33}--\eqref{7.36} generalize  Buschman--Erdelyi ones \eqref{71}--\eqref{72} respectively. They reduce to Buschman--Erdelyi for the choice of parameters $a=-(\nu+\mu),~ b= 1+ \nu - \mu,~ c=1-\mu$. For operators \eqref{7.33}--\eqref{7.36} the above results are generalized with necessary changes. For example they are factorized via more simple operators \eqref{7.37}--\eqref{7.38} with special choice of parameters.

Operators \eqref{7.37}--\eqref{7.38} are generalizations of  \eqref{73}--\eqref{733}. For them the next result is true.

\theor Operators  \eqref{7.37}--\eqref{7.38} may be extended to isometric in $L_2 (0, \infty)$ if and only if they coincide with   Buschman--Erdelyi operators of zero order smoothness \eqref{73}--\eqref{733} for integer values of $\nu=\frac{1}{2}(b-a-1)$.

This theorem single out Buschman--Erdelyi operators of zero order smoothness at least in the class
\eqref{7.33}--\eqref{7.38}. Operators \eqref{7.33}--\eqref{7.36} are generalizations of fractional integrals. Analogically may be studied generalizations to \eqref{6.1}--\eqref{6.2}, \eqref{6.6}, \eqref{6.16}--\eqref{6.17}.

More general are operators with $G$ function kernel.

\begin{eqnarray}
& & {{_1}G_{0+}}(\alpha, \beta, \delta, \gamma)[f]= \frac{2^{\delta}}{\Gamma(1-\alpha)\Gamma(1-\beta)}\cdot \label{7.39} \\
& & \int\limits_0^x (\frac{x}{y}-1)^{-\delta}(\frac{x}{y}+1)^{1+\delta-\alpha-\beta} G_{2~2}^{1~2} \lr{\frac{x}{2y}-\frac{1}{2}|_{\gamma, \, \delta}^{\alpha, \, \beta}}  f(y)\, dy. \nonumber
\end{eqnarray}
Another operators are with different interval of integration and parameters of $G$ function. For $\alpha=1-a,~ \beta=1-b, \delta=1-c, \gamma=0$ \eqref{7.39} reduce to \eqref{7.33}, for $\alpha=1+\nu,~ \beta=-\nu, \delta= \gamma=0$  \eqref{7.39} reduce to Buschman--Erdelyi operators of zero order smoothness
${{_1}S_{0+}^{\nu}}$.

Further generalizations are in terms of Wright or Fox functions. They lead to Wright--Buschman--Erdelyi and Fox--Buschman--Erdelyi operators. These classes are connected with Sonine--Dimovski and Poisson--Dimovski transmutations  \cite{Dim},  \cite{Kir}, and also with generalized fractional integrals introduced by Kiryakova \cite{Kir}.

\bigskip

\subsection{Application of Buschman--Erdelyi transmutations in works of V.\,Katrakhov.}
\medskip

V.\,Katrakhov found a new approach for boundary value problems for elliptic equations with strong singularities of infinite order. For example for Poisson equation he studied problems with solutions of arbitrary growth. At singular point he proposed the new kind of boundary condition: $K$---trace. His results are based on constant usage of Buschman--Erdelyi transmutations of the first kind  for definition of norms, solution estimates and correctness proofs \cite{Kat1}--\cite{Kat2}.

Moreover in joint papers with I.\,Kipriyanov he introduced and studied new classes of pseudodifferential operators based on transmutational technics \cite{Kat3}. These results were
paraphrased in reorganized manner also in \cite{Car2}.

\bigskip

\addcontentsline{toc}{section}{\ \ \ \ \  References\hspace{104mm}  28}

\newpage

--------------------------------------------------------------------------------------------

\selectlanguage{russian}

\bigskip
--------------------------------------------------------------------------------------------
\newpage

\section*{\\ОПЕРАТОРЫ ПРЕОБРАЗОВАНИЯ БУШМАНА--ЭРДЕЙИ,\\
ИСТОРИЧЕСКИЙ СВЕДЕНИЯ, КЛАССИФИКАЦИЯ, ПРИЛОЖЕНИЯ.}

\section*{\\Русская версия статьи.}

\section*{\\С.М. СИТНИК
\\Воронежский институт МВД, Воронеж, Россия.}

\addcontentsline{toc}{section}{РУССКАЯ ВЕРСИЯ СТАТЬИ.}

\section{\\Введение.}

\subsection{Операторы преобразования (ОП).}

Теория операторов преобразования---это существенное обобщение теории подобия конечномерных матриц. Дадим сразу основное определение.

Definition 1. {Пусть дана пара операторов $(A,B)$. Оператор
$T$ называется \textit{оператором преобразования} (ОП, сплетающий оператор, transmutation, intertwining operator), если
на элементах подходящих функциональных пространств выполняется соотношение
\begin{equation}
\label{1.1}
\Large {T\,A=B\,T.}
\end{equation}
}

Ясно, что понятие ОП является прямым и далеко идущим обобщением понятия подобия матриц из линейной алгебры. Но ОП \textit{не сводятся к подобным (или эквивалентным) операторам},  так как сплетаемые операторы как правило являются неограниченными в естественных пространствах, к тому же обратный к ОП не обязан существовать, действовать в том же пространстве или быть ограниченным. Так что спектры операторов, сплетаемых ОП, как правило не совпадают. Кроме того, сами ОП могут быть неограниченными. Так бывает, например, в теории преобразований Дарбу, предметом которой является нахождение дифференциальных операторов преобразования (подстановок или замен) между парой дифференциальных же операторов, таким образом в этом случае все три рассматриваемых оператора являются неограниченными в естественных пространствах. При этом теория преобразований Дарбу как соответствующий раздел теории дифференциальных уравнений также вписывается в общую схему теории операторов преобразования при её расширенном понимании. Кроме того, можно рассматривать операторы преобразования не только для пары дифференциальных операторов. В теории ОП встречаются задачи для  следующих разнообразных типов операторов: интегральных,  интегро--дифференциальных, дифференциально--разностных (например,  типа Дункла), дифференциальных или интегро--дифференциальных   бесконечного порядка (например, в вопросах, связанных с леммой Шура о дополняемости), общих линейных  в фиксированных функциональных пространствах, псевдодифференциальных и операторно--дифференциальных (абстрактных дифференциальных).

Возможность, чтобы исходная и преобразованная функции принадлежали различным пространствам, что принято подчёркивать использованием различных обозначений для переменных, позволяет включить в общую схему ОП все классические интегральные преобразования: Фурье, Лапласа (на самом деле Петцваля), Меллина, Ханкеля, Вейерштрасса, Конторовича--Лебедева, Мелера--Фока, Станковича и другие. В общую схему ОП также включаются конечные интегральные преобразования Г.А.\,Гринберга.

В квантовой физике при рассмотрении уравнения Шрёдингера и задач теории рассеяния встречается специальный класс ОП --- волновые операторы.

Коммутирующие операторы любой природы также подходят под определение ОП. Наиболее близко к духу и задачам теории ОП относится изучение операторов, коммутирующих с производными. Сами ОП в этом случае зачастую представляются формальными рядами, псевдо--дифференциальными операторами или дифференциальными операторами бесконечного порядка. Описание коммутантов напрямую связано с описанием всего семейства ОП для заданной пары по его единственному представителю. В этом классе задач фундаментальные приложения нашла теория операторных свёрток, особенно свертки Берга--Димовски.
Начинают находить приложения в теории ОП и результаты для коммутирующих дифференциальных операторов, восходящие к классическим работам Бёрчнела и Чонди (J.L.\,Burchnall, T.W.\,Chaundy). Теория ОП также связана с вопросами факторизации дифференциальных операторов.

Отдельный класс ОП составляют преобразования, которые для одного и того же уравнения связывают краевые условия различных типов, например, Неймана и Дирихле.

Как же обычно используются операторы преобразования? Пусть, например,  мы изучаем некоторый достаточно сложно устроенный оператор $A$. При этом нужные свойства уже известны для модельного более простого оператора $B$. Тогда, если существует ОП (\ref{1.1}), то часто удаётся перенести свойства модельного оператора  $B$ и на $A$. Такова в нескольких словах примерная схема типичного использования ОП в конкретных задачах.

В частности, если рассматривается уравнение $Au=f$ с оператором $A$, то применяя к нему ОП $T$ со сплетающим свойством (\ref{1.1}), получаем уравнение с оператором $B$ вида $Bv=g$, где обозначено $v=Tu, g=Tf$. Поэтому, если второе уравнение с оператором $B$ является более простым, и для него уже известны формулы для решений, то мы получаем и представления для решений первого уравнения $u=T^{-1}v$. Разумеется, при этом обратный оператор преобразования должен существовать и действовать в рассматриваемых пространствах, а для получения явных представлений решений должно быть получено и явное представление этого обратного оператора. Таково одно из простейших применений техники ОП в теории дифференциальных уравнений, как обыкновенных, так и  с частными производными.

Изложению теории ОП и их приложениям полностью посвящены монографии  \cite{Car1}-\cite{Tri1}, а также подробный обзор автора \cite{Sit1}. Методы теории ОП составляют кроме того существенные части монографий \cite{CarSho}-\cite{Hro}, сейчас этот список можно дополнить почти до 100 монографий.

Следует специально отметить  монографию Д.К.\,Фаге и Н.И.\,Нагнибида \cite{FaNa}. В этой монографии практически никак не отражены уже известные к тому времени результаты теории ОП, что полностью компенсируется изложением в основном собственных результатов авторов по одной из самых трудных задач теории ОП --- их построении  для дифференциальных операторов высоких порядков с переменными коэффициентами. Кроме того, в эту монографию вошли и многие другие вопросы: решение задачи об операторах, коммутирующих с производными в пространствах аналитических функций (включая исправление ошибочных результатов Дельсарта и Лионса), создание законченной теории разрешимости для уравнения Бианки, теория операторно--аналитических функций (первоначально возникшая в работах В.А.\,Марченко), исследование операторов дифференцирования, интегрирования и корней из них в пространствах аналитических функций.

Сделаем одно терминологическое замечание. В западной литературе принят для ОП термин \textit{"transmutation"}, восходящий к Ж.\,Дельсарту. Как отмечает Р.\,Кэрролл, похожий термин "transformation"\  при этом закрепляется за классическими интегральными преобразованиями Фурье, Лапласа, Меллина, Ханкеля и другими подобными им.  Приведём дословную цитату из \cite{Car3}: "Such operators are often called transformation operators by the Russian school (Levitan, Naimark, Marchenko et. al.), but transformation seems too broad a term, and, since some of the machinery seems "magical"\  at times, we have followed Lions and Delsarte in using the word transmutation".

 В настоящее время теория операторов преобразования представляет собой полностью оформившийся самостоятельный раздел математики, находящийся на стыке дифференциальных и интегральных уравнений, функционального анализа, теории функций, комплексного анализа, теории специальных функций и дробного интегродифференицрования.
 Необходимость теории операторов преобразования  доказана большим числом её приложений. Методы операторов преобразования применяются в теории обратных задач, определяя обобщённое преобразование Фурье, спектральную функцию и решения знаменитого уравнения Левитана; в теории рассеяния через операторы преобразования выписывается не менее знаменитое уравнение Марченко; в  спектральной теории получаются известные формулы следов и асимптотика спектральной функции; оценки ядер операторов преобразования отвечают за устойчивость обратных задач и задач рассеяния; в теории нелинейных дифференциальных уравнений метод Лакса использует операторы преобразования для доказательства существования решений и построения солитонов. Определёнными разновидностями операторов преобразования являются части теорий  обобщённых аналитических функций, операторов обобщённого сдвига и обобщённых операторных свёрток, метод преобразования Дарбу. В теории уравнений с частными производными методы операторов преобразования применяются для построения явных выражений для решений возмущённых  задач через решения невозмущённых, изучении сингулярных и вырождающихся краевых задач, псевдодифференциальных операторов,  задач для решений с существенными особенностями на части границы  во внутренних или угловых точках, оценки скорости убывания решений некоторых эллиптических и ультраэллиптических уравнений. Теория операторов преобразования позволяет дать новую классификацию специальных функций и интегральных операторов со специальными функциями в ядрах, в том числе различных операторов дробного интегродифференцирования. В теории функций найдены приложения операторов преобразования к вложениям функциональных пространств и обобщению операторов Харди, расширению теории Пэли--Винера, построению различных конструкций обобщённого сдвига и основанным на них обобщённых вариантов гармонического анализа.  Методы теории операторов преобразования с успехом применяются во многих прикладных задачах: оценках решений Йоста в квантовой теории рассеяния, обратных задачах, исследовании системы Дирака и других матричных систем дифференциальных уравнений, операторных и дифференциально--операторных уравнениях, различных интегральных уравнениях, в том числе со специальными функциями в ядрах,  теории вероятностей и случайных процессов, линейном стохастическом оценивании, фильтрации, стохастических случайных уравнениях, обратных задачах геофизики и трансзвуковой газодинамики. Кроме уже известных для метода Лакса и преобразований Дарбу всё время увеличивается число новых приложений ОП к нелинейным дифференциальным уравнениям и исследованию солитонов.

 Фактически современная теория операторов преобразования возникла из двух примеров, ставших классическими \cite{Sit1}.
Первым примером являются ОП, переводящие оператор Штурма--Лиувилля с некоторым потенциалом $q(x)$ во вторую производную:
\begin{equation}
T(D^2\,y(x)+q(x)y(x))=D^2\,(Ty(x)), D^2\,y(x)=y''(x),
\end{equation}
при некотором выборе естественных краевых условий.

Второй пример---это  задача о преобразовании оператора Бесселя во вторую производную:
\begin{equation}
T\lr{B_\nu} f=\lr{D^2} Tf, B_{\nu}=D^2+\frac{2\nu +1}{x}D, D^2=\frac{d^2}{dx^2}, \nu \in \mathbb{C}.
\end{equation}
На этом пути возникли ОП Сонина--Пуассона--Дельсарта, Бушмана--Эрдейи и их многочисленные обобщения.
Такие операторы преобразования находят многочисленные приложения при изучении одного класса уравнений с частными производными с особенностями, типичным представителем которого является $B$--эллиптическое уравнение с операторами Бесселя по каждой переменной вида
\begin{equation}
\label{59}
\sum_{k=1}^{n}B_{\nu,x_k}u(x_1,\dots, x_n)=f,
\end{equation}
аналогично рассматриваются $B$--гиперболические и $B$--параболические уравнения. Изучение этого класса уравнений было начато в работах Эйлера, Пуассона, Дарбу, продолжено в теории обобщённого осесимметрического потенциала А.\,Вайнштейна и в трудах отечественных математиков И.Е.\,Егорова, Я.И.\,Житомирского, Л.Д.\,Кудрявцева, П.И.\,Лизоркина, М.И.\,Матийчука, Л.Г.\,Михайлова, М.Н.\,Олевского, М.М.\,Смирнова, С.А.\,Терсенова,   Хе Кан Чера, А.И.\,Янушаускаса   и других.

Наиболее полно весь круг вопросов для  уравнений с операторами Бесселя был изучен воронежским математиком И.А.\,Киприяновым и его учениками Л.А.\,Ивановым, А.В.\,Рыжковым, В.В.\,Катраховым, В.П.\,Архиповым, А.Н.\,Байдаковым, Б.М.\,Богачёвым, А.Л.\,Бродским, Г.А.\,Виноградовой, В.А.\,Зайцевым, Ю.В.\,Засориным, Г.М.\,Каганом,  А.А.\,Катраховой, Н.И.\,Киприяновой, В.И.\,Кононенко, М.И.\,Ключанцевым, А.А.\,Куликовым, А.А.\,Лариным, М.А.\,Лейзиным, Л.Н.\,Ляховым, А.Б.\,Муравником, И.П.\,Половинкиным, А.Ю.\,Сазоновым, С.М.\,Ситником, В.П.\,Шацким, В.Я.\,Ярославцевой;   основные результаты этого направления представлены в \cite{Kip1}. Для описания классов решений соответствующих уравнений И.А.\,Киприяновым были введены и изучены функциональные пространства \cite{Kip2},  позднее названные его именем.

В  этом направлении работал В.В.\,Катрахов, сейчас  уравнения с оператором Бесселя и связанные с ними вопросы изучают  А.В.\,Глушак, В.С.\,Гулиев,  Л.Н.\,Ляхов   со своими коллегами и учениками. Задачи для операторно--дифференциальных (абстрактных) уравнений вида (\ref{59}), берущие начало в известной монографии \cite{CarSho}, рассматривали А.В.\,Глушак, С.Б.\,Шмулевич, В.Д.\,Репников и другие.

Операторы преобразования являются одним из основных инструментов изучения этого класса уравнений. Они используются для получения формул для решений, фундаментальных решений, описания особенностей, постановок краевых задач и в других вопросах.

Кратко опишем структуру предлагаемой статьи. Данная работа носит в основном обзорный характер. Вместе с тем основной результат об условиях ограниченности и унитарности операторов Бушмана--Эрдейи приведён с полным доказательством (теорема \ref{tnorm}). Во вводном первом разделе излагаются  исторические и приоритетные сведения. Используется предложенная автором удобная классификация различных классов операторов Бушмана--Эрдейи. На основе этой классификации во втором разделе изложены основные результаты автора по операторам преобразования Бушмана--Эрдейи первого рода, включая операторы нулевого порядка гладкости, в  третьем разделе---по операторам преобразования Бушмана--Эрдейи второго рода, в четвёртом---по операторам преобразования Бушмана--Эрдейи третьего рода, а также унитарным операторам преобразования Сонина--Катрахова и Пуассона--Катрахова. В заключительном пятом разделе приведены приложения операторов
преобразования Бушмана--Эрдейи различных классов к вложению пространств И.А.Киприянова в весовые пространства С.Л.Соболева, формулам для решений уравнений с частными производными с операторами Бесселя, уравнениям Эйлера--Пуассона--Дарбу, включая лемму Копсона, построению операторов обобщённого сдвига, операторам Дункла, Преобразованию Радона, построению обобщённых сферических гармоник и $B$--гармонических полиномов, а также доказательству унитарности в пространстве Лебега обобщений классических операторов Харди. В заключение этого раздела кратко перечислены результаты В.В.Катрахова по приложению операторов преобразования  Бушмана--Эрдейи к построение нового класса псевдодифференциальных операторов и изучению введённого им класса краевых задач с $K$---следом с существенными особенностями в решениях.

Сделаем одно техническое замечание. Каждый конкретный результат о том, что некоторый ОП сплетает пару преобразуемых операторов, должен сопровождаться аккуратным указанием всех областей определения операторов и классов рассматриваемых функций. Для краткости эта конкретика в данной статье указывается не для каждого результата, хотя для каждой теоремы такие точные условия получены. Кроме того, мы используем в некоторых случаях для краткости термин "оператор"\  вместо более точного "дифференциальное выражение".

\subsection{Операторы преобразования Бушмана--Эрдейи.}

Название  \textsl{"операторы Бушмана--Эрдейи"} \ было предложено автором, в последнее время оно стало общепринятым. Интегральные уравнения с подобными операторами рассматривались с середины 1950--х годов, автором было впервые показано, что операторы Бушмана--Эрдейи являются операторами преобразования для дифференциального выражения Бесселя и изучены их специальные свойства именно как операторов преобразования.  Частными случаями операторов преобразования Бушмана--Эрдейи являются известные классические операторы преобразования Сонина и Пуассона, а их обобщениями являются операторы преобразования Сонина--Димовски и Пуассона--Димовски для гипербесселевых уравнений и функций.

Операторы Бушмана--Эрдейи имеют многочисленные модификации. Автором предложена удобная классификация их различных вариантов. Операторы Бушмана--Эрдейи первого рода содержат ядра, выражающиеся через функции Лежандра первого рода. Их предельным случаем являются операторы нулевого порядка гладкости, играющие важную роль в различных приложениях. Операторы Бушмана--Эрдейи второго рода содержат ядра, выражающиеся через функции Лежандра второго рода. Комбинация операторов первого и второго родов приводит к операторам Бушмана--Эрдейи третьего рода. При специальном выборе параметров они сводятся к унитарным операторам преобразованиям, которые автор назвал унитарными операторами преобразования Сонина--Катрахова и Пуассона--Катрахова, в честь  В.В.\,Катрахова, начавшего их изучение.

Изучение разрешимости и обратимости данных операторов было начато в 1960--х годах в работах    Р.\,Бушмана и А.\,Эрдейи \cite{Bu1}--\cite{Er2}. Операторы Бушмана--Эрдейи или их аналоги изучались также в работах T.P.\,Higgins, Ta Li, E.R.\,Love, G.M.\,Habibullah, K.N.\,Srivastava, Динь Хоанг Ань, В.И.\,Смирнова, Н.А.\,Вирченко, И.\,Федотовой, А.А.\,Килбаса, О.В.\,Скоромник и др. При этом  изучались задачи о решении интегральных уравнений с этими операторами, их факторизации и обращения. Основные результаты указанного периода изложены в монографии \cite{KK}, хотя случай выбранных нами пределов интегрирования считается там особым и не рассматривается, за исключением одного набора формул композиции, некоторые результаты для особого выбора пределов были добавлены в английское расширенное издание монографии \cite{KK}.

Наиболее полное изучение операторов Бушмана--Эрдейи на наш взгляд было проведено в работах автора в 1980--1990-е годы \cite{Sit3}--\cite{Sit4}, и затем продолжено в \cite{Sit2}--\cite{Sit10} и ряде других работ. При этом необходимо отметить, что  роль операторов Бушмана--Эрдейи как ОП до работ  \cite{Sit3}--\cite{Sit4} вообще ранее нигде не отмечалась и не рассматривалась.

 Из работ, в которых изучались операторы Бушмана--Эрдейи как интегральные операторы, отметим работы Н.А.\,Вирченко, А.А.\,Килбаса  и их учеников. Так в работах А.А.\,Килбаса и О.В.\,Скоромник \cite{KiSk1}--\cite{KiSk2} рассматривается действие операторов Бушмана--Эрдейи в весовых пространствах Лебега, а также  многомерные обобщения в виде интегралов по пирамидальным областям. В монографии Н.А.\,Вирченко и И.\,Федотовой \cite{ViFe} вводятся некоторые обобщения стандартных функций Лежандра, а затем рассматриваются напоминающие операторы Бушмана--Эрдейи, но не содержащие их как частные случаи, интегральные операторы с введёнными функциями в ядрах на всей положительной полуоси (операторы Бушмана--Эрдейи определены на части положительной полуоси).

Важность операторов Бушмана--Эрдейи во многом обусловлена их многочисленными приложениями. Например, они встречаются в следующих вопросах теории уравнений с частными производными: при решении задачи Дирихле для уравнения Эйлера--Пуассона--Дарбу в четверти плоскости и установлении соотношений между значениями решений уравнения Эйлера--Пуассона--Дарбу на многообразии начальных данных и характеристике, теории преобразования Радона, так как в силу результатов Людвига  действие преобразования Радона при разложении по сферическим гармоникам сводится как раз к операторам  Бушмана--Эрдейи по радиальной переменной, при исследовании краевых задач для различных уравнений с существенными особенностями внутри области, доказательству вложения пространств И.А.\,Киприянова  в весовые пространства С.Л.\,Соболева, установлению связей между операторами преобразования и волновыми операторами теории рассеяния, обобщению классических интегральных представлений Сонина и Пуассона и операторов преобразования Сонина--Пуассона--Дельсарта.

\section{\\Операторы преобразования Бушмана--Эрдейи первого рода.}

\setcounter{equation}{3}

\subsection{Операторы преобразования Сонина--Пуассона--Дельсарта.}

Рассмотрим, наверное, самый известный класс ОП, сплетающих дифференциальный оператор Бесселя со второй производной:
\begin{equation}
\label{2.1}
T\lr{B_\nu} f=\lr{D^2} Tf, B_{\nu}=D^2+\frac{2\nu +1}{x}D, D^2=\frac{d^2}{dx^2}, \nu \in \mathbb{C}.
\end{equation}

Definition 2.
 ОП Пуассона называется выражение
\begin{equation}
\label{2.2}
P_{\nu}f=\frac{1}{\Gamma(\nu+1)2^{\nu}x^{2\nu}}
\int_0^x \left( x^2-t^2\right)^{\nu-\frac{1}{2}}f(t)\,dt,\Re \nu> -\frac{1}{2}.
\end{equation}
ОП Сонина называется выражение
\begin{equation}
\label{2.3}
S_{\nu}f=\frac{2^{\nu+\frac{1}{2}}}{\Gamma(\frac{1}{2}-\nu)}\frac{d}{dx}
\int_0^x \left( x^2-t^2\right)^{-\nu-\frac{1}{2}}t^{2\nu+1}f(t)\,dt,\Re \nu< \frac{1}{2}.
\end{equation}
Операторы (\ref{2.2})--(\ref{2.3}) действуют как ОП по формулам
\begin{equation}
\label{56}
S_\nu B_\nu=D^2 S_\nu, P_\nu D^2=B_\nu P_\nu.
\end{equation}
Их можно доопределить на все значения $\nu\in\mathbb{C}$.
Исторически более точно называть введённые операторы именами Сонина--Пуассона--Дельсарта. Большинство математиков узнали об этих операторах из  статьи Б.М. Левитана \cite{Lev3}.

Важным обобщением операторов Сонина--Пуассона--Дельсарта являются ОП для гипербесселевых функций. Теория таких функций была первоначально заложена в работах Куммера и Делерю. Полное исследование гипербесселевых функций, дифференциальных уравнений для них и соответствующих операторов преобразования  было исчерпывающе проведено в работах И.\,Димовски и его учеников \cite{Dim}. Соответствующие ОП получили в литературе названия ОП Сонина--Димовски и Пуассона--Димовски, они также изучались в работах ученицы И.\,Димовски --- В.\,Киряковой \cite{Kir}. В теории гипербесселевых функций, дифференциальных уравнений и операторов преобразования для них центральную роль играет знаменитое интегральное преобразование Обрешкова, введённое болгарским математиком Н.\,Обрешковым \cite{Kir}. Это преобразование, ядро которого выражается в общем случае через $G$--функцию Майера, является одновременным обобщением преобразований Лапласа, Меллина, синус-- и косинус преобразований Фурье, Ханкеля, Майера и других классических интегральных преобразований. Различные формы гипербесселевых функций, дифференциальных уравнений и  операторов преобразований для них, а также частные случаи преобразования Обрешкова многократно впоследствии переоткрывались, этот процесс продолжается и до настоящего времени. По мнению автора, преобразование Обрешкова, наряду с преобразованиями Фурье, Меллина, Лапласа, Станковича относится к небольшому числу фундаментальных преобразований Анализа, из которых, как из кирпичиков, складываются многие другие преобразования, а также основанные на них конструкции и приложения.

\subsection{Определения и основные свойства.}

Теперь перейдём к описанию основных свойств операторов преобразования Бушмана--Эрдейи.
Это класс ОП, который при определённом выборе параметров является одновременным обобщением ОП СПД и их сопряжённых, операторов дробного интегродифференцирования Римана--Лиувилля
и Эрдейи--Кобера, а также интегральных преобразований Мелера--Фока.

Definition 3.
Операторами Бушмана--Эрдейи первого рода называются интегральные операторы
\begin{eqnarray}
\label{71}
B_{0+}^{\nu,\mu}f=\int_0^x \left( x^2-t^2\right)^{-\frac{\mu}{2}}P_\nu^\mu \left(\frac{x}{t}\right)f(t)d\,t,\\
E_{0+}^{\nu,\mu}f=\int_0^x \left( x^2-t^2\right)^{-\frac{\mu}{2}}\mathbb{P}_\nu^\mu \left(\frac{t}{x}\right)f(t)d\,t,\\
B_{-}^{\nu,\mu}f=\int_x^\infty \left( t^2-x^2\right)^{-\frac{\mu}{2}}P_\nu^\mu \left(\frac{t}{x}\right)f(t)d\,t,\\
\label{72}
E_{-}^{\nu,\mu}f=\int_x^\infty \left( t^2-x^2\right)^{-\frac{\mu}{2}}\mathbb{P}_\nu^\mu \left(\frac{x}{t}\right)f(t)d\,t.\\ \nonumber
\end{eqnarray}
Здесь $P_\nu^\mu(z)$---функция Лежандра первого рода, $\mathbb{P}_\nu^\mu(z)$---та же функция на разрезе $-1\leq  t \leq 1$, $f(x)$---локально суммируемая функция, удовлетворяющая некоторым ограничениям на рост при $x\to 0,x\to\infty$. Параметры $\mu,\nu$---комплексные числа, $\Re \mu <1$, можно ограничиться значениями $\Re \nu \geq -1/2$.

Приведём основные результаты, в основном следуя в изложении \cite{Sit3}, \cite{Sit4}, а также  \cite{Sit1}, \cite{Sit2}. Все рассмотрения ведутся ниже на полуоси. Поэтому будем обозначать через $L_2$ пространство $L_2(0, \infty)$ и $L_{2, k}$ весовое пространство $L_{2, k}(0, \infty)$ со степенным весом и нормой
\begin{equation}
\int_0^\infty |f(x)|^2 x^{2k+1}\,dx.
\end{equation}
$\mathbb{N}$ \ обозначает множество натуральных, $\mathbb{N}_0$--неотрицательных целых, $\mathbb{Z}$--целых и $\mathbb{R}$--действительных чисел.

Вначале распространим определение  3  на   случай значения параметра $\mu =1$.

Definition 4.
 Введём при $\mu =1$ операторы  Бушмана--Эрдейи нулевого порядка гладкости по формулам
\begin{eqnarray}
\label{73}
B_{0+}^{\nu,1}f=\frac{d}{dx}\int_0^x P_\nu \left(\frac{x}{t}\right)f(t)\,dt,\\
\label{731}
E_{0+}^{\nu,1}f=\int_0^x P_\nu \left(\frac{t}{x}\right)\frac{df(t)}{dt}\,dt,\\
\label{732}
B_{-}^{\nu,1}f=\int_x^\infty P_\nu \left(\frac{t}{x}\right)(-\frac{df(t)}{dt})\,dt,\\
\label{733}
E_{-}^{\nu,1}f=(-\frac{d}{dx})\int_x^\infty P_\nu \left(\frac{x}{t}\right)f(t)\,dt,
\end{eqnarray}
где $P_\nu(z)=P_\nu^0(z)$---функция Лежандра.

\theor{\it
 Справедливы следующие формулы факторизации операторов Бушмана--Эрдейи на подходящих функциях через дробные интегралы Римана--Лиувилля и Бушмана--Эрдейи нулевого порядка гладкости:
\begin{equation}\label{1.9}
{B_{0+}^{\nu,\,\mu} f=I_{0+}^{1-\mu}~ {_1 S^{\nu}_{0+}f},~B_{-}^{\nu, \,\mu} f={_1 P^{\nu}_{-}}~ I_{-}^{1-\mu}f,}
\end{equation}
\begin{equation}\label{1.10}
{E_{0+}^{\nu,\,\mu} f={_1 P^{\nu}_{0+}}~I_{0+}^{1-\mu}f,~E_{-}^{\nu, \, \mu} f= I_{-}^{1-\mu}~{_1 S^{\nu}_{-}}f.}
\end{equation}
}

Эти формулы позволяют "разделить"\  параметры $\nu$ и $\mu$. Мы докажем, что операторы
\eqref{73}--\eqref{733} являются изоморфизмами пространств $L_2(0, \infty)$, если $\nu$
не равно некоторым исключительным значениям. Поэтому операторы \eqref{71}--\eqref{72}
по действию в пространствах типа $L_2$ в определённом смысле подобны операторам дробного интегродиффенцирования  $I^{1-\mu}$, с которыми они совпадают при $\nu=0$. Далее операторы Бушмана--Эрдейи будут доопределены при всех значениях $\mu$.

Исходя из этого, введём следующее

Definition 5.{ Число $\rho=1-Re\,\mu $ назовём порядком гладкости операторов Бушмана--Эрдейи \eqref{71}--\eqref{72}.}

Таким образом, при $\rho > 0$ (то есть при $Re\, \mu > 1$) операторы Бушмана--Эрдейи
являются сглаживающими, а при $\rho < 0$ (то есть при $Re\, \mu < 1$) уменьшающими
гладкость в пространствах типа $L_2 (0, \infty)$. Операторы \eqref{73}--\eqref{733},
для которых $\rho = 0$, являются по данному определению операторами нулевого порядка гладкости.

Перечислим основные свойства операторов Бушмана--Эрдейи первого рода
\eqref{71}--\eqref{72} с функцией Лежандра I рода в ядре.
При некоторых специальных значениях параметров $\nu,~\mu$ операторы Бу\-шмана--Эрдейи
сводятся к более простым. Так при значениях $\mu=-\nu$ или $\mu=\nu+2$ они являются
операторами Эрдейи--Кобера; при $\nu = 0$ операторами дробного интегродифференцирования
$I_{0+}^{1-\mu}$ или $I_{-}^{1-\mu}$; при $\nu=-\frac{1}{2}$, $\mu=0$ или $\mu=1$
ядра выражаются через эллиптические интегралы; при  $\mu=0$,  $x=1$, $v=it-\frac{1}{2}$  оператор $B_{-}^{\nu, \, 0}$ лишь на постоянную отличается от преобразования Мелера--Фока.

Будем рассматривать наряду с оператором Бесселя также тесно связанный с ним дифференциальный оператор
\begin{equation}
\label{75}
L_{\nu}=D^2-\frac{\nu(\nu+1)}{x^2}=\left(\frac{d}{dx}-\frac{\nu}{x}\right)
\left(\frac{d}{dx}+\frac{\nu}{x}\right),
\end{equation}
который при $\nu \in \mathbb{N}$ является оператором углового момента из квантовой механики.
Их взаимосвязь устанавливает

\theor \label{tOP}
{\it Пусть пара ОП $X_\nu, Y_\nu$ сплетают $L_{\nu}$ и вторую  производную:
\begin{equation}
\label{76}
X_\nu L_{\nu}=D^2 X_\nu , Y_\nu D^2 = L_{\nu} Y_\nu.
\end{equation}
Введём новую пару ОП по формулам
\begin{equation}
\label{77}
S_\nu=X_{\nu-1/2} x^{\nu+1/2}, P_\nu=x^{-(\nu+1/2)} Y_{\nu-1/2}.
\end{equation}
Тогда пара новых ОП $S_\nu, P_\nu$ сплетают оператор Бесселя и вторую производную:
\begin{equation}
\label{78}
S_\nu B_\nu = D^2 S_\nu, P_\nu D^2 = B_\nu P_\nu.
\end{equation}
}

\theor\label{t6}
{\it Пусть $Re \, \mu \leq 1$. Тогда оператор $B_{0+}^{\nu, \, \mu}$ является оператором преобразования типа Сонина и удовлетворяет на подходящих функциях соотношению \eqref{76}.
}

Аналогичный результат справедлив и для других операторов Бушмана--Эрдейи. При этом
$E_{-}^{\nu, \, \mu}$ также является оператором типа Сонина, а $E_{0+}^{\nu, \, \mu}$
и $B_{-}^{\nu, \, \mu}$ -- операторами типа Пуассона.

Теперь сделаем важное замечание. Из приведённой теоремы следует, что ОП Бушмана--Эрдейи связывают собственные функции операторов Бесселя и второй производной. Таким образом, половина ОП Бушмана--Эрдейи переводят тригонометрические или экспоненциальные функции в приведённые функции Бесселя, а другая половина наоборот. Эти формулы здесь не приводятся, их нетрудно выписать явно. Все они являются обобщениями исходных формул Сонина и Пуассона   и представляют существенный интерес. Ещё раз отметим, что подобные формулы являются непосредственными следствиями доказанных сплетающих свойств ОП Бушмана--Эрдейи, и могут быть непосредственно проверены при помощи таблиц интегралов от специальных функций.

Перейдём к вопросу о различных факторизациях операторов Бушмана--Эрдейи через операторы
Эрдейи--Кобера и дробные интегралы Римана--Лиувилля.

Приведём список основных операторов дробного интегродифференцирования: Римана--Лиувилля, Эрдейи--Кобера, дробного интеграла по произвольной функции $g(x)$, см. \cite{KK}
\begin{eqnarray}
\label{61}
I_{0+,x}^{\alpha}f=\frac{1}{\Gamma(\alpha)}\int_0^x \left( x-t\right)^{\alpha-1}f(t)d\,t,\\ \nonumber
I_{-,x}^{\alpha}f=\frac{1}{\Gamma(\alpha)}\int_x^\infty \left( t-x\right)^{\alpha-1}f(t)d\,t,
\end{eqnarray}
\begin{eqnarray}
\label{62}
I_{0+,2,\eta}^{\alpha}f=\frac{2 x^{-2\lr{\alpha+\eta}}}{\Gamma(\alpha)}\int_0^x \left( x^2-t^2\right)^{\alpha-1}t^{2\eta+1}f(t)d\,t,\\ \nonumber
I_{-,2,\eta}^{\alpha}f=\frac{2 x^{2\eta}}{\Gamma(\alpha)}\int_x^\infty \left( t^2-x^2\right)^{\alpha-1}t^{1-2\lr{\alpha+\eta}}f(t)d\,t,
\end{eqnarray}
\begin{eqnarray}
\label{63}
I_{0+,g}^{\alpha}f=\frac{1}{\Gamma(\alpha)}\int_0^x \left( g(x)-g(t)\right)^{\alpha-1}g'(t)f(t)d\,t,\\ \nonumber
I_{-,g}^{\alpha}f=\frac{1}{\Gamma(\alpha)}\int_x^\infty \left( g(t)-g(x)\right)^{\alpha-1}g'(t)f(t)d\,t,
\end{eqnarray}
во всех случаях предполагается, что $\Re\alpha>0$, на оставшиеся значения
$\alpha$ формулы также без труда продолжаются \cite{KK}.
При этом обычные дробные интегралы получаются при выборе в (\ref{63}) $g(x)=x$, Эрдейи--Кобера при $g(x)=x^2$, Адамара при $g(x)=\ln x$.

\theor\label{tfact7} Справедливы следующие формулы факторизации операторов
Бушмана--Эрдейи первого рода через операторы дробного интегродифференцирования и Эрдейи--Кобера:
\begin{eqnarray}
& & B_{0+}^{\nu, \, \mu}=I_{0+}^{\nu+1-\mu} I_{0+; \, 2, \, \nu+ \frac{1}{2}}^{-(\nu+1)} {\lr{\frac{2}{x}}}^{\nu+1}\label{2.17}, \\
& & E_{0+}^{\nu, \, \mu}= {\lr{\frac{x}{2}}}^{\nu+1} I_{0+; \, 2, \, - \frac{1}{2}}^{\nu+1} I_{0+}^{-(\nu+\mu)}  \label{2.18}, \\
& & B_{-}^{\nu, \, \mu}= {\lr{\frac{2}{x}}}^{\nu+1}I_{-; \, 2, \, \nu+ 1}^{-(\nu+1)} I_{-}^{\nu - \mu+2}  \label{2.19}, \\
& & E_{-}^{\nu, \, \mu}= I_{-}^{-(\nu+\mu)} I_{-; \, 2, \, 0} ^{\nu+1} {\lr{\frac{x}{2}}}^{\nu+1}  \label{2.20}.
\end{eqnarray}

Частными случаями введённых операторов являются и определённые выше операторы преобразования Сонина--Пуассона--Дельсарта.

Теперь рассмотрим более подробно свойства ОП Бушмана--Эрдейи нулевого порядка гладкости, введённых по формулам (\ref{73}). Подобный им оператор был построен В.В. Катраховым  путём домножения стандартного ОП Сонина на обычный дробный интеграл с целью взаимно компенсировать гладкость этих двух операторов и получить новый, который бы действовал в одном пространстве типа $L_2(0,\infty)$.

Напомним \cite{Mar}, что преобразованием Меллина функции $f(x)$ называется функция $g(s)$, которая определяется по формуле
\begin{equation}
\label{710}
g(s)=M{f}(s)=\int_0^\infty x^{s-1} f(x)\,dx.
\end{equation}
Определим также свёртку Меллина
\begin{equation}
\label{711}
(f_1*f_2)(x)=\int_0^\infty  f_1\left(\frac{x}{y}\right) f_2(y)\,\frac{dy}{y},
\end{equation}
при этом оператор свёртки с ядром $K$ действует в образах преобразования Меллина как умножение на мультипликатор
\begin{eqnarray}
\label{con}
M{Af}(s)=\int_0^\infty  K\left(\frac{x}{y}\right) f(y)\,\frac{dy}{y}=M{K*f}(s)
=m_A(s)M{f}(s),\\\nonumber m_A(s)=M{K}(s).\phantom{1111111111111111111}
\end{eqnarray}

Заметим, что преобразование Меллина является обобщённым преобразованием Фурье на полуоси по мере Хаара  $\frac{dy}{y}$ \cite{Hel}. Его роль велика в теории специальных функций, например, гамма--функция является преобразованием Меллина экспоненты. С преобразованием Меллина связан важный прорыв в 1970--х годах, когда в основном усилиями О.И.\,Маричева была полностью доказана и приспособлена для нужд вычисления интегралов известная теорема Д.Л.\,Слейтер, позволяющая для большинства образов преобразований Меллина восстановить оригинал в явном виде по простому алгоритму через гипергеометрические функции \cite{Mar}.  Эта теорема стала основой универсального мощного метода вычисления интегралов, который  позволил решить многие задачи в теории дифференциальных и интегральных уравнений, а также воплотился в передовые технологии символьного интегрирования пакета MATHEMATICA фирмы Wolfram Research.

\theor\label{tnorm}
Оператор Бушмана--Эрдейи нулевого порядка гладкости $B_{0+}^{\nu,1}$, определённый по формуле (\ref{73}), действует в образах преобразования Меллина как свёртка \ref{con} с мультипликатором
\begin{equation}\label{712}
m(s)=\frac{\Gamma(-s/2+\frac{\nu}{2}+1)\Gamma(-s/2-\frac{\nu}{2}+1/2)}
{\Gamma(1/2-\frac{s}{2})\Gamma(1-\frac{s}{2})}\label{s1}
\end{equation}
при условии $\Re s<\min(2+\Re\nu,1-\Re\nu)$. Для его нормы, которая является периодической функцией по $\nu$, справедлива формула
\begin{equation}\label{s2}
\|B_{0+}^{\nu,1}\|_{L_2}=\frac{1}{\min(1,\sqrt{1-\sin\pi\nu})}.
\end{equation}
Оператор ограничен в $L_2(0,\infty)$ при $\nu\neq 2k+1/2, k\in \mathbb{Z}$
и неограничен при выполнении условия $\nu= 2k+1/2, k\in \mathbb{Z}$.

\vspace{3mm}

Эта теорема является основной в работе, поэтому дадим её полное доказательство.

1. Вначале докажем формулу (\ref{712}) с нужным мультипликатором. Используя последовательно формулы (7), с.~130, (2) с.~129, (4) с.~130 из \cite{Mar}, получим
$$
M(B_{0+}^{\nu,1})(s)=\frac{\Gamma(2-s)}{\Gamma(1-s)}M
\left[\int_0^\infty
\left\{H(\frac{x}{y}-1)P_\nu (\frac{x}{y})
\right\}
\left\{y f(y)\right\}\frac{dy}{y}
\right](s-1)=
$$
$$
=\frac{\Gamma(2-s)}{\Gamma(1-s)} M
\left[(x^2-1)_+^0P_\nu^0 (x)
\right]
(s-1)M\left[f\right](s),
$$
где использованы обозначения из \cite{Mar} для функции Хевисайда и усечённой степенной функции
$$
x_+^\alpha=\left\{
\begin{array}{rl}
x^\alpha, & \mbox{если } x\geqslant 0 \\
0, & \mbox{если } x<0 \\
\end{array}\right.
,\ H(x)=x_+^0=\left\{
\begin{array}{rl}
1, & \mbox{если } x\geqslant 0 \\
0, & \mbox{если } x<0. \\
\end{array}\right.
$$
Далее, используя формулы 14(1) с.~234 и 4 с.~130 из \cite{Mar}, получаем
$$
M\left[(x-1)_+^0 P_\nu^0 (\sqrt x)
\right](s)=
\frac{\Gamma(\frac{1}{2}+\frac{\nu}{2}-s)\Gamma(-\frac{\nu}{2}-s)}
{\Gamma(1-s)\Gamma(\frac{1}{2}-s)},
$$
$$
M\left[(x^2-1)_+^0 P_\nu^0 (x)
\right](s-1)=\frac{1}{2}\cdot
\frac
{
\Gamma(\frac{1}{2}+\frac{\nu}{2}-\frac{s-1}{2})
\Gamma(-\frac{\nu}{2}-\frac{s-1}{2})
}
{
\Gamma(1-\frac{s-1}{2})\Gamma(\frac{1}{2}-\frac{s-1}{2})
}=
$$
$$
=\frac{1}{2}\cdot\frac
{
\Gamma(-\frac{s}{2}+\frac{\nu}{2}+1)
\Gamma(-\frac{s}{2}-\frac{\nu}{2}+\frac{1}{2})
}
{\Gamma(-\frac{s}{2}+\frac{3}{2})\Gamma(-\frac{s}{2}+1)}
$$
при условиях $\Re s<\min(2+\Re\nu,1-\Re\nu)$. Отсюда выводим формулу для мультипликатора
$$
M(B_{0+}^{\nu,1})(s)=\frac{1}{2}\cdot\frac{\Gamma(2-s)}{\Gamma(1-s)}\cdot
{\Gamma(-\frac{s}{2}+\frac{3}{2})\Gamma(-\frac{s}{2}+1)}.
$$
Применяя к $\Gamma(2-s)$ формулу Лежандра удвоения аргумента гамма--функции (см., например, \cite{BE1}), получим
$$
M(B_{0+}^{\nu,1})(s)=\frac{2^{-s}}{\sqrt\pi}\cdot
\frac
{
\Gamma(-\frac{s}{2}+\frac{\nu}{2}+1)
\Gamma(-\frac{s}{2}-\frac{\nu}{2}+\frac{1}{2})
}
{\Gamma(1-s)}.
$$
Ещё одно применение формулы удвоения Лежандра к $\Gamma(1-s)$ приводит к нужной формуле для мультипликатора (\ref{s1}).

В работе \cite{Sit3} показано, что за счёт рассмотрения подходящих факторизаций, условия справедливости доказанной формулы, которые являются завышенными, так как выводились сразу для всего класса обобщённых гипергеометрических функций Гаусса, можно несколько расширить. В частности, в нужном нам конкретном случае ядра с функцией Лежандра формула для мультипликатора справедлива при условиях $0<\Re s<1$ при всех значениях параметра $\nu$.

2. Теперь установим формулу для нормы (\ref{s2}). Из найденной формулы для мультипликатора в силу теоремы 4.7 из \cite{Sit1} получаем на прямой $\Re s=1/2, s=i u+1/2$
$$
|M(B_{0+}^{\nu,1})(i u+1/2)|=\frac{1}{\sqrt{2\pi}}\left|\frac
{
\Gamma(-i\frac{u}{2}-\frac{\nu}{2}+\frac{1}{4})
\Gamma(-i\frac{u}{2}+\frac{\nu}{2}+\frac{3}{4})
}
{\Gamma(\frac{1}{2}-iu)}\right|.
$$
Далее будем опускать у мультипликатора указание на порождающий его оператор. Используем формулу для модуля комплексного числа $|z|=\sqrt{z\bar{z}}$ и тождество для гамма--функции $\overline{\Gamma(z)}=\Gamma(\bar z)$, вытекающее из её определения в виде интеграла. Последнее равенство справедливо для класса так называемых вещественно--аналитических функций, к которому относится и гамма--функция. Тогда получим
$$
|M(B_{0+}^{\nu,1})(i u+1/2)|=
$$
$$
=\frac{1}{\sqrt{2\pi}}\left|\frac
{
\Gamma(-i\frac{u}{2}-\frac{\nu}{2}+\frac{1}{4})
\Gamma(i\frac{u}{2}-\frac{\nu}{2}+\frac{1}{4})
\Gamma(-i\frac{u}{2}+\frac{\nu}{2}+\frac{3}{4})
\Gamma(i\frac{u}{2}+\frac{\nu}{2}+\frac{3}{4})
}
{\Gamma(\frac{1}{2}-iu)\Gamma(\frac{1}{2}+iu)}\right|.
$$
В числителе объединим крайние и средние сомножители, и три образовавшиеся пары гамма--функций преобразуем по известной формуле (см. \cite{BE1})
$$
\Gamma(\frac{1}{2}+z)\  \Gamma(\frac{1}{2}-z)=\frac{\pi}{\cos \pi z}.
$$
В результате получим
$$
|M(B_{0+}^{\nu,1})(i u+1/2)|=
\sqrt{
\frac{\cos(\pi i u)}
{2\cos\pi(\frac{\nu}{2}+\frac{1}{4}+i\frac{u}{2})
\cos\pi(\frac{\nu}{2}+\frac{1}{4}-i\frac{u}{2})}
}=
$$
$$
=\sqrt{
\frac{\ch(\pi i u)}{\ch\pi u-\sin\pi\nu}
}
$$
Далее обозначим $t=\ch\pi u, 1\leqslant t <\infty$. Отсюда, применяя опять условие из теоремы 4.7 \cite{Sit1}, получаем
$$
\sup_{u\in\mathbb{R}} |m(i u+\frac{1}{2})|=\sup_{1\leqslant t <\infty}
\sqrt{
\frac{t}{t-\sin\pi\nu}
}.
$$
Поэтому, если $\sin\pi\nu\geqslant 0$, то супремум достигается при $t=1$, и справедлива формула (\ref{s2}) для нормы
$$
\|B_{0+}^{\nu,1}\|_{L_2}=\frac{1}{\sqrt{1-\sin\pi\nu}}.
$$
Если же $\sin\pi\nu\leqslant 0$, то супремум достигается при $t\to\infty$,
и справедлива формула
$$
\|B_{0+}^{\nu,1}\|_{L_2}=1.
$$
Эта часть теоремы доказана.

3. Условия ограниченности или неограниченности следуют из найденной формулы для нормы и условий процитированной теоремы 4.7. Периодичность нормы по параметру $\nu$ очевидна из найденного явного выражения для нормы.

Теорема полностью доказана.

Приведём формулы для мультипликаторов всех операторов Бушмана--Эрдейи нулевого порядка гладкости.

\theor\label{tmult}
Операторы Бушмана--Эрдейи нулевого порядка гладкости действуют в образах преобразования Меллина как свёртки по формуле (\ref{con}). Для их мультипликаторов справедливы формулы:

\begin{eqnarray}
& & m_{{_1S_{0+}^{\nu}}}(s)=\frac{\Gamma(-\frac{s}{2}+\frac{\nu}{2}+1) \Gamma(-\frac{s}{2}-\frac{\nu}{2}+\frac{1}{2})}{\Gamma(\frac{1}{2}-\frac{s}{2})\Gamma(1-\frac{s}{2})}= \label{3.11}; \\
& & =\frac{2^{-s}}{\sqrt{ \pi}} \frac{\Gamma(-\frac{s}{2}-\frac{\nu}{2}+\frac{1}{2}) \Gamma(-\frac{s}{2}+\frac{\nu}{2}+1)}{\Gamma(1-s)} , Re\, s < \min (2 + Re \, \nu, 1- Re\, \nu) \nonumber ; \\
& & m_{{_1P_{0+}^{\nu}}}(s)=\frac{\Gamma(\frac{1}{2}-\frac{s}{2})\Gamma(1-\frac{s}{2})}{\Gamma(-\frac{s}{2}+\frac{\nu}{2}+1) \Gamma(-\frac{s}{2}-\frac{\nu}{2}+\frac{1}{2})},~ Re\, s < 1; \label{3.12} \\
& & m_{{_1P_{-}^{\nu}}}(s)=\frac{\Gamma(\frac{s}{2}+\frac{\nu}{2}+1) \Gamma(\frac{s}{2}-\frac{\nu}{2})}{\Gamma(\frac{s}{2})\Gamma(\frac{s}{2}+\frac{1}{2})}, Re \, s > \max(Re \, \nu, -1-Re\, \nu); \label{3.13} \\
& & m_{{_1S_{-}^{\nu}}}(s)=\frac{\Gamma(\frac{s}{2})\Gamma(\frac{s}{2}+\frac{1}{2})}{\Gamma(\frac{s}{2}+\frac{\nu}{2}+\frac{1}{2}) \Gamma(\frac{s}{2}-\frac{\nu}{2})}, Re \, s >0 \label{3.14}
\end{eqnarray}

Справедливы следующие формулы для норм операторов Бушмана--Эрдейи нулевого порядка гладкости в $L_2$:
\begin{eqnarray}
& & \| _1{S_{0+}^{\nu}} \| = \| _1{P_{-}^{\nu}}\|= 1/ \min(1, \sqrt{1- \sin \pi \nu}), \label{3.22} \\
& & \| _1{P_{0+}^{\nu}}\| = \| _1{S_{-}^{\nu}}\|= \max(1, \sqrt{1- \sin \pi \nu}). \label{3.23}
\end{eqnarray}

Аналогичные результаты получены в \cite{Sit2}--\cite{Sit3} и для пространств со степенным весом.

Следствие 1. Нормы операторов \eqref{73} -- \eqref{733} периодичны по $\nu$ с периодом 2, то есть $\|X^{\nu}\|=\|X^{\nu+2}\|$, где $X^{\nu}$ --- любой из операторов \eqref{73} -- \eqref{733}.

Следствие 2. Нормы операторов ${_1 S_{0+}^{\nu}}$, ${_1 P_{-}^{\nu}}$ не ограничены в совокупности по $\nu$, каждая из этих норм не меньше 1. Если $\sin \pi \nu \leq 0$, то эти нормы равны 1. Указанные операторы неограничены в $L_2$ тогда и только тогда, когда $\sin \pi \nu = 1$ (или $\nu=(2k) + 1/2,~k \in \mathbb{Z}$).

Следствие 3. Нормы операторов ${_1 P_{0+}^{\nu}}$, ${_1 S_{-}^{\nu}}$
ограничены в совокупности по $\nu$, каждая из этих норм не больше $\sqrt{2}$. Все эти операторы ограничены в $L_2$ при всех $\nu$. Если $\sin \pi \nu \geq 0$, то их $L_2$ -- норма равна 1. Максимальное значение нормы, равное $\sqrt 2 $, достигается тогда и только тогда, когда $\sin \pi \nu = -1$ (или $\nu= -1/2+(2k) ,~k \in \mathbb{Z}$).

Важнейшим свойством операторов Бушмана--Эрдейи нулевого порядка гладкости является их унитарность при целых $\nu$. Отметим, что при интерпретации $L_{\nu}$ как оператора углового момента в квантовой механике, параметр $\nu$ как раз и принимает целые неотрицательные значения.

\theor Для унитарности в $L_2$ операторов \eqref{73} -- \eqref{733} необходимо и достаточно, чтобы число $\nu$ было целым. В этом случае пары операторов
$({_1 S_{0+}^{\nu}}$, ${_1 P_{-}^{\nu}})$ и  $({_1 S_{-}^{\nu}}$, ${_1 P_{0+}^{\nu}})$
взаимно обратны.

Перед формулировкой одного частного случая  предположим, что операторы \eqref{73} -- \eqref{733} заданы на таких функциях $f(x)$, что возможно в определениях выполнить внешнее дифференцирование или под знаком интеграла интегрирование по частям (для этого достаточно предположить, что $x f(x) \to 0$ при $x \to 0$). Тогда при $\nu=1$
\begin{equation}\label {3.25}
{_1{P_{0+}^{1}}f=(I-H_1)f,~_1{S_{-}^{1}}f=(I-H_2)f,}
\end{equation}
где $H_1,~ H_2$ -- операторы Харди,
\begin{equation}\label {3.26}
{H_1 f = \frac{1}{x} \int\limits_0^x f(y) dy,~H_2 f = \int\limits_x^{\infty}  \frac{f(y)}{y} dy,}
\end{equation}
$I$ --- единичный оператор.

Следствие 4. Операторы \eqref{3.25} являются унитарными взаимно обратными в $L_2$ операторами. Они сплетают дифференциальные выражения $d^2 / d x^2$ и $d^2 / d x^2 - 2/ x^2$.

Унитарность в $L_2$ сдвинутых операторов Харди \eqref{3.25} известна, см. \cite{KuPe}. Ниже в приложениях описаны новые подобные унитарные операторы.

Далее перечислим некоторые общие свойства операторов, которые  действуют по формуле \eqref{con} как умножение на некоторый мультипликатор в образах преобразования  Меллина и одновременно являются сплетающими для второй производной и оператора углового момента.

\theor\label{tOPmult} Пусть оператор $S_{\nu}$ действует по формулам \eqref{con}
и \eqref{76}. Тогда

а) его мультипликатор удовлетворяет функциональному уравнению

\begin{equation}\label{5.1}{m(s)=m(s-2)\frac{(s-1)(s-2)}{(s-1)(s-2)-\nu(\nu+1)};}
\end{equation}

б) если функция $p(s)$ периодична с периодом 2 (то есть $p(s)=p(s-2)$),  то функция $p(s)m(s)$ является мультипликатором нового оператора  преобразования $S_2^{\nu}$, опять же сплетающего $L_{\nu}$ и вторую производную по правилу
\eqref{76}.

Последняя теорема ещё раз показывает, насколько полезно изучение ОП в терминах мультипликаторов преобразования Меллина.

Определим преобразование Стилтьеса (см., например, \cite{KK}) по формуле
$$
(S f)(x)= \int\limits_0^{\infty} \frac{f(t)}{x+t} dt.
$$
Этот оператор также действует по формуле \eqref{con} с мультипликатором  $p(s)= \pi /sin (\pi s)$ и ограничен в $L_2$. Очевидно, что $p(s)=p(s-2)$. Поэтому из теоремы \ref{tOPmult} следует, что композиция преобразования Стилтьеса с ограниченными сплетающими операторами \eqref{73}--\eqref{733} снова является оператором преобразованием
того же типа, ограниченным в $L_2$.

На этом получены многочисленные новые явные семейства ОП, ядра которых выражаются через различные специальные функции.

\section{\\Операторы преобразования Бушмана--Эрдейи второго рода.}

\setcounter{equation}{41}

Теперь определим и изучим операторы Бушмана--Эрдейи второго рода.

Definition 5.
Введём новую пару операторов Бушмана--Эрдейи с функциями Лежандра второго рода  в ядре
\begin{equation}\label{6.1}
{{_2S^{\nu}}f=\frac{2}{\pi} \left( - \int\limits_0^x (x^2-y^2)^{-\frac{1}{2}}Q_{\nu}^1 (\frac{x}{y}) f(y) dy  +
\int\limits_x^{\infty} (y^2-x^2)^{-\frac{1}{2}}\mathbb{Q}_{\nu}^1 (\frac{x}{y}) f(y) dy\right),}
\end{equation}

\begin{equation}\label{6.2} {{_2P^{\nu}}f=\frac{2}{\pi} \left( - \int\limits_0^x (x^2-y^2)^{-\frac{1}{2}}\mathbb{Q}_{\nu}^1 (\frac{y}{x}) f(y) dy  -
\int\limits_x^{\infty} (y^2-x^2)^{-\frac{1}{2}}Q_{\nu}^1 (\frac{y}{x}) f(y) dy\right).}
\end{equation}

Эти операторы являются аналогами операторов первого рода нулевого порядка гладкости.

При $y \to x \pm 0$ интегралы понимаются в смысле главного значения. Отметим без доказательства, что эти операторы определены и являются сплетающими при некоторых условиях на функции $f(x)$  (при этом оператор \eqref{6.1} будет типа Сонина, \eqref{6.2}  --- типа Пуассона).

\theor  Операторы  \eqref{6.1} -- \eqref{6.2} представимы в виде \eqref{con}
с мультипликаторами
\begin{eqnarray}
& & m_{_2S^{\nu}}(s)=p(s) \ m_{_1S_{-}^{\nu}}(s), \label{6.3}\\
& & m_{_2P^{\nu}}(s)=\frac{1}{p(s)} \ m_{_1P_{-}^{\nu}}(s), \label{6.4}
\end{eqnarray}
где мультипликаторы операторов ${_1S_-^{\nu}}$, ${_1P_-^{\nu}}$ определены формулами \eqref{3.13} -- \eqref{3.14}, а функция $p(s)$ (с периодом 2) равна

\begin{equation}\label {6.5}{p(s)=\frac{\sin \pi \nu+ \cos \pi s}{\sin \pi \nu - \sin \pi s}.}
\end{equation}

\theor Справедливы формулы для норм
\begin{eqnarray}
& &  \| {_2S^{\nu}} \|_{L_2}= \max (1, \sqrt {1+\sin \pi \nu}) , \label{6.9} \\
& &  \| {_2P^{\nu}} \|_{L_2}= 1 / {\min (1, \sqrt {1+\sin \pi \nu})} . \label{6.10}
\end{eqnarray}

Следствие.  Оператор ${_2S^{\nu}}$ ограничен при всех $\nu$. Оператор ${_2P^{\nu}}$ не является непрерывным тогда и только тогда, когда
$\sin \pi \nu=-1$.

\theor Для унитарности в $L_2$ операторов ${_2S^{\nu}}$ и ${_2P^{\nu}}$ необходимо и достаточно, чтобы параметр $\nu$ был целым
числом.

\theor Пусть $\nu=i \beta+1/2,~\beta \in \mathbb{R}$. Тогда
\begin{equation}\label{6.11}{\| {_2S^{\nu}} \|_{L_2}=\sqrt{1+\ch \pi \beta},~\| {_2P^{\nu}} \|_{L_2}=1.}
\end{equation}

\theor Справедливы представления
\begin{eqnarray}
& & {_2S^0} f = \frac{2}{\pi} \int\limits_0^{\infty} \frac{y}{x^2-y^2}f(y)\,dy, \label{6.12} \\
& & {_2S^{-1}} f = \frac{2}{\pi} \int\limits_0^{\infty} \frac{x}{x^2-y^2}f(y)\,dy. \label{6.13}
\end{eqnarray}

Таким образом, в этом случае оператор ${_2S^{\nu}}$ сводится к паре известных преобразований Гильберта на полуоси \cite{KK}.

Для операторов второго рода введём также более общие аналоги операторов преобразования Бушмана--Эрдейи первого рода с двумя параметрами по формуле:

\begin{eqnarray}
& & {_2S^{\nu,\mu}}f=\frac{2}{\pi} \left( \int\limits_0^x (x^2+y^2)^{-\frac{\mu}{2}} e^{-\mu \pi i} Q_{\nu}^{\mu}( \frac{x}{y}) f(y)\, dy\right. + \label{6.6} \\
& & +\left. \int\limits_x^{\infty} (y^2+x^2)^{-\frac{\mu}{2}}\mathbb{Q}_{\nu}^{\mu} (\frac{x}{y}) f(y)\, dy\right) \nonumber,
\end{eqnarray}
где $Q_{\nu}^{\mu}(z)$ --- функция Лежандра второго рода, $\mathbb{Q}_{\nu}^{\mu}(z)$ --- значение этой функции на разрезе, $Re\, \nu < 1$. Второй подобный оператор определяется как формально сопряжённый в $L_2(0,\infty)$ к (\ref{6.6}).

\theor
На функциях из $C_0^{\infty}(0, \infty)$ оператор \eqref{6.6} определён и действует по формуле
$$
M{[_2S^{\nu}]}(s)=m(s)\cdot M[x^{1-\mu} f](s), \label{6.7}
$$
\begin{eqnarray}
 m(s)=2^{\mu-1} \left( \frac{ \cos \pi(\mu-s) - \cos \pi \nu}{ \sin \pi(\mu-s) - \sin \pi \nu}  \right) \cdot\\
\cdot \left( \frac{\Gamma(\frac{s}{2})\Gamma(\frac{s}{2}+\frac{1}{2}))}{\Gamma(\frac{s}{2}+\frac{1-\nu-\mu}{2}) \Gamma(\frac{s}{2}+1+\frac{\nu-\mu}{2})} \right). \nonumber
\end{eqnarray}

\section{\\Операторы преобразования Бушмана--Эрдейи третьего рода.}

\setcounter{equation}{53}

\subsection{Операторы преобразования Сонина--Катрахова и Пуассона--Катрахова.}

Перейдём к построению операторов преобразования, унитарных при  всех $\nu$. Такие операторы определяются по формулам:
\begin{eqnarray}
& & S_U^{\nu} f = - \sin \frac{\pi \nu}{2}\  {_2S^{\nu}}f+ \cos \frac{\pi \nu}{2}\  {_1S_-^{\nu}}f, \label{6.14} \\
& & P_U^{\nu} f = - \sin \frac{\pi \nu}{2}\  {_2P^{\nu}}f+ \cos \frac{\pi \nu}{2}\  {_1P_-^{\nu}}f. \label{6.15}
\end{eqnarray}
Для любых значений $\nu \in \mathbb{R}$ они являются линейными комбинациями операторов преобразования Бушмана--Эрдейи 1 и 2 рода нулевого порядка гладкости. Их можно отнести к операторам Бушмана--Эрдейи третьего рода (см. ниже). В интегральной форме эти операторы имеют вид:

\begin{eqnarray}
& & S_U^{\nu} f = \cos \frac{\pi \nu}{2} \left(- \frac{d}{dx} \right) \int\limits_x^{\infty} P_{\nu}\lr{\frac{x}{y}} f(y)\,dy +  \label{6.16}\\
& & + \frac{2}{\pi} \sin \frac{\pi \nu}{2} \left(  \int\limits_0^x (x^2-y^2)^{-\frac{1}{2}}Q_{\nu}^1 \lr{\frac{x}{y}} f(y)\,dy  \right.-
%\nonumber\\& & -
 \int\limits_x^{\infty} (y^2-x^2)^{-\frac{1}{2}}\mathbb{Q}_{\nu}^1 \lr{\frac{x}{y}} f(y)\,dy \Biggl. \Biggr), \nonumber \\
& & P_U^{\nu} f = \cos \frac{\pi \nu}{2}  \int\limits_0^{x} P_{\nu}\lr{\frac{y}{x}} \left( \frac{d}{dy} \right) f(y)\,dy - \label{6.17} \\
& &  -\frac{2}{\pi} \sin \frac{\pi \nu}{2} \left( - \int\limits_0^x (x^2-y^2)^{-\frac{1}{2}}\mathbb{Q}_{\nu}^1\lr{\frac{y}{x}} f(y)\,dy   \right.-
%\nonumber\\& & -
\int\limits_x^{\infty} (y^2-x^2)^{-\frac{1}{2}} Q_{\nu}^1 \lr{\frac{y}{x}} f(y)\,dy \Biggl. \Biggr). \nonumber
\end{eqnarray}

\theor Операторы \eqref{6.14}--\eqref{6.15},  \eqref{6.16}--\eqref{6.17} при всех $\nu\in\mathbb{R}$ являются унитарными, взаимно сопряжёнными и обратными в $L_2$. Они являются
сплетающими и действуют по формулам \eqref{75}. При этом $S_U^{\nu}$ является оператором типа Сонина (Сонина--Катрахова), а $P_U^{\nu}$ --- типа Пуассона (Пуассона--Катрахова).

ОП в форме подобной (\ref{6.16})--(\ref{6.17}), но только с ядрами, выражающимися через общую гипергеометрическую функцию Гаусса, были впервые построены в 1980 г. В.В.\,Катраховым. Поэтому автор предлагает названия: операторы преобразования Сонина--Катрахова и Пуассона--Катрахова. Их выражение через функции Лежандра первого и второго родов получено автором, кроме того их удаётся включить в общую схему построения операторов преобразования композиционным методом \cite{SiKa2}, \cite{SiKa3}, \cite{Sit5}. При этом основными становятся наиболее простые  формулы факторизации вида \eqref{6.14}--\eqref{6.15}. На этом пути построение подобных операторов перестаёт быть специальным искусным приёмом, а встраивается в общую методику построения целых классов подобных операторов преобразования композиционным методом.

\subsection{Операторы преобразования Бушмана--Эрдейи третьего рода с произвольной весовой функцией.}

Рассмотрим синус и косинус-преобразования Фурье и обратные к ним

\begin{equation}\label{1}{F_c f = \sqrt{\frac{2}{\pi}} \int\limits_0^{\infty} f(y) \cos (t y) \, dy,
\ \ F_c^{-1}=F_c,}
\end{equation}

\begin{equation}\label{2}{F_s f = \sqrt{\frac{2}{\pi}} \int\limits_0^{\infty} f(y) \sin (t y) \, dy,}
\ \ F_s^{-1}=F_s.
\end{equation}

Определим преобразование Фурье-Бесселя по формулам

\begin{eqnarray}
& &  F_{\nu} f = \frac{1}{2^{\nu} \Gamma (\nu+1 )} \int\limits_0^{\infty} f(y)\, j_{\nu}(t y) \, y^{2 \nu + 1} \, dy = \nonumber\\
& & = \int\limits_0^{\infty} f(y) \frac{J_\nu(t y)}{(t y)^{\nu}}\,  y^{2 \nu + 1} \, dy = \frac{1}{t^{\nu}} \int\limits_0^{\infty} f(y)  J_{\nu}(t y) \, y^{\nu + 1} \, dy, \label{3} \\
& & F_{\nu}^{-1} f = \frac{1}{(y)^{\nu}} \int\limits_0^{\infty} f(t) J_{\nu}(y t)\, t^{\nu + 1} \, dt. \label{4}
\end{eqnarray}
Здесь $J_\nu(\cdot)$---обычная \cite{BE1}, а\  $j_\nu(\cdot)$---нормированная \cite{Kip1} функции Бесселя.
Операторы \eqref{1}-\eqref{2} самосопряженные унитарные в $L_2(0, \infty)$.  Операторы \eqref{3}-\eqref{4} самосопряженные унитарные в $L_{2,\, \nu}(0, \infty)$.

Определим первую пару операторов преобразования Бушмана--Эрдейи третьего рода на подходящих функциях с произвольной весовой функцией по формулам
\begin{eqnarray}
& &  S_{\nu,\, c}^{(\varphi)} = F^{-1}_c \left( \frac{1}{\varphi (t)} F_{\nu} \right), \label{5} \\
& &  P_{\nu,\, c}^{(\varphi)} = F^{-1}_{\nu} \left( \varphi (t) F_{c} \right), \label{6}
\end{eqnarray}

и вторую пару по формулам

\begin{eqnarray}
& &  S_{\nu,\, s}^{(\varphi)} = F^{-1}_s \left( \frac{1}{\varphi (t)} F_{\nu} \right), \label{7} \\
& &  P_{\nu,\, s}^{(\varphi)} = F^{-1}_{\nu} \left( \varphi (t) F_{s} \right). \label{8},
\end{eqnarray}
где $\varphi (t)$---произвольная весовая функция.

Введённые операторы преобразования на подходящих функциях сплетают $B_{\nu}$ и $D^2$, можно дать их интегральное представление.

\theor Определим операторы преобразования, сплетающие $B_{\nu}$ и $D^2$, по формулам

$$
S^{(\varphi)}_{\nu, \left\{ \begin{matrix} s \\ c \end{matrix} \right\} }=
F^{-1}_{ \left\{ \begin{matrix} s \\ c \end{matrix} \right\} } \left( \frac{1}{\varphi(t)}F_{\nu} \right),
$$

$$
P^{(\varphi)}_{\nu, \left\{ \begin{matrix} s \\ c \end{matrix} \right\} }=
F^{-1}_{ \nu } \left( \varphi(t)F_{ \left\{ \begin{matrix} s \\ c \end{matrix} \right\} } \right).
$$

Тогда для оператора типа Сонина справедливо представление (формальное)

\begin{equation}{\left( S^{(\varphi)}_{\nu, \left\{ \begin{matrix} s \\ c \end{matrix} \right\} } f \right) (x) = \sqrt{\frac{2}{\pi}} \int\limits_0^{\infty} K(x, y) f(y) \, dy,}
\end{equation}
где
$$
K(x, y) = y^{\nu + 1} \int\limits_0^{\infty} \frac{\left\{ \begin{matrix} \sin (x t) \\ \cos (x t) \end{matrix} \right\}}{\varphi (t)\, t^{\nu} } J_{\nu} (y t) dt.
$$

Представление для оператора типа Пуассона имеет вид

\begin{equation}{\left( P^{(\varphi)}_{\nu, \left\{ \begin{matrix} s \\ c \end{matrix} \right\} } f \right) (x) = \sqrt{\frac{2}{\pi}} \int\limits_0^{\infty} G(x, y) f(y) \, dy,}
\end{equation}
где
$$
G(x, y) = \frac{1}{x^\nu} \int\limits_0^{\infty} \varphi (t) \, t^{\nu + 1} \left\{ \begin{matrix} \sin (y t) \\ \cos (y t) \end{matrix}\right\}  J_{\nu} (x t) dt.
$$

Частным случаем введённых операторов третьего рода являются определённые выше унитарные операторы Сонина--Катрахова и Пуассона--Катрахова, которые получаются при выборе весовой функции $\varphi (t)$ в виде некоторой зависящей от параметра $\nu$ степени. Достаточно полная теория и приложения операторов преобразования Бушмана--Эрдейи третьего рода планируются автором к публикации в отдельной статье.

\section{\\Некоторые приложения операторов преобразования Бушмана--Эрдейи.}

\setcounter{equation}{57}

В этом разделе приведём приложения введённых операторов. По необходимости мы ограничиваемся только формулировками основных результатов или даже только простым перечислением фактов с минимальным набором ссылок, некоторые приложения только намечены, их исследование пока подробно не проводилось.

\bigskip

\subsection{Оценки норм в пространствах И.А. Киприянова.}

Рассмотрим множество функций $\mathbb{D}(0, \infty)$. Если $f(x) \in \mathbb{D}(0, \infty)$, то $f(x) \in C^{\infty}(0, \infty),~ f(x)$ ---финитна на бесконечности. На этом множестве функций введём полунормы
\begin{eqnarray}
& & \|f\|_{h_2^{\alpha}}=\|D_-^{\alpha}f\|_{L_2(0, \infty)} \label{7.1} \\
& &  \|f\|_{\widehat{h}_2^{\alpha}}=\|x^{\alpha} (-\frac{1}{x}\frac{d}{dx})^{\alpha}f\|_{L_2(0, \infty)} \label{7.2}
\end{eqnarray}
где $D_-^{\alpha}$ --- дробная производная Римана--Лиувилля, оператор в \eqref{7.2} определяется по формуле
\begin{equation}\label{7.3}{(-\frac{1}{x}\frac{d}{dx})^{\beta}=2^{\beta}I_{-; \, 2, \,0}^{-\beta}x^{-2 \beta},}
\end{equation}
$I_{-; 2, \, 0}^{-\beta}$ --- оператор Эрдейи--Кобера, $\alpha$ --- произвольное действительное число. При $\beta = n \in \mathbb{N}_0$ выражение \eqref{7.3} понимается в обычном смысле, что согласуется с определением выше.

\theor
 Пусть $f(x) \in \mathbb{D}(0, \infty)$. Тогда справедливы тождества:
\begin{eqnarray}
& & D_-^{\alpha}f={_1S_-^{\alpha-1}} {x^{\alpha} (-\frac{1}{x}\frac{d}{dx})^{\alpha}} f, \label{7.4} \\
& & x^{\alpha} (-\frac{1}{x}\frac{d}{dx})^{\alpha}f={_1P_-^{\alpha-1}} D_-^{\alpha}f. \label{7.5}
\end{eqnarray}

Таким образом, операторы Бушмана--Эрдейи нулевого порядка гладкости первого рода осуществляют связь между дифференциальными операторами (при $\alpha \in \mathbb{N}$) из определений полунорм \eqref{7.1} и \eqref{7.2}.

\theor
 Пусть $f(x) \in \mathbb{D}(0, \infty)$. Тогда справедливы неравенства между полунормами
\begin{eqnarray}
& &  \|f\|_{h_2^{\alpha}} \leq \max (1, \sqrt{1+\sin \pi \alpha}) \|f\|_{\widehat{h}_2^{\alpha}}, \label{7.7}\\
& & \|f\|_{\widehat{h}_2^{\alpha}} \leq \frac{1}{\min (1, \sqrt{1+\sin \pi \alpha})} \|f\|_{h_2^{\alpha}}, \label{7.8}
\end{eqnarray}
где $\alpha$ --- любое действительное число, $\alpha \neq -\frac{1}{2}+2k,~k \in \mathbb{Z}$.

Постоянные в неравенствах \eqref{7.7}--\eqref{7.8} не меньше единицы, что будет далее использовано. В случае
$\sin \pi \alpha = -1 $ или $\alpha = -\frac{1}{2}+2k,~k \in \mathbb{Z}$, оценка \eqref{7.8} не имеет места.

Введём на $\mathbb{D} (0, \infty )$ соболевскую норму
\begin{equation}\label{7.9}{\|f\|_{W_2^{\alpha}}=\|f\|_{L_2 (0, \infty)}+\|f\|_{h_2^{\alpha}}.}
\end{equation}
Введём также другую норму
\begin{equation}\label{7.10}{\|f\|_{\widehat{W}_2^{\alpha}}=\|f\|_{L_2 (0, \infty)}+\|f\|_{\widehat{h}_2^{\alpha}}}
\end{equation}
Пространства $W_2^{\alpha},~ \widehat{W}_2^{\alpha}$ определим как замыкания $D(0,
\infty)$ по нормам \eqref{7.9} и \eqref{7.10} соответственно.

\theor а) при всех $\alpha \in \mathbb{R}$ пространство $\widehat{W}_2^{\alpha}$ непрерывно вложено в $W_2^{\alpha}$, причём
\begin{equation}\label{7.11}{\|f\|_{W_2^{\alpha}}\leq A_1 \|f\|_{\widehat{W}_2^{\alpha}},}
\end{equation}
где $A_1=\max (1, \sqrt{1+\sin \pi \alpha})$.

б) Пусть $\sin \pi \alpha \neq -1$ или $\alpha \neq -\frac{1}{2} + 2k, ~ k \in \mathbb{Z}.  $ Тогда справедливо обратное вложение $W_2^{\alpha}$  в $\widehat{W}_2^{\alpha}$, причём
\begin{equation}\label{7.12}{\|f\|_{\widehat{W}_2^{\alpha}}\leq A_2 \|f\|_{W_2^{\alpha}},}
\end{equation}
где $A_2 =1/  \min (1, \sqrt{1+\sin \pi \alpha})$.

в) Пусть $\sin \pi \alpha \neq -1$, тогда пространства $W_2^{\alpha}$  и $\widehat{W}_2^{\alpha}$ изоморфны, а их нормы эквивалентны.

г) Константы в неравенствах вложений \eqref{7.11}--\eqref{7.12} точные.

Эта теорема фактически является следствием результатов по ограниченности операторов Бушмана--Эрдейи нулевого порядка гладкости в $L_2$. В свою очередь, из теоремы об унитарности этих операторов вытекает

\theor Нормы
\begin{eqnarray}
& & \|f\|_{W_2^{\alpha}} = \sum\limits_{j=0}^s \| D_-^j f\|_{L_2}, \label{7.13} \\
& & \|f\|_{\widehat{W}_2^{\alpha}}=\sum\limits_{j=0}^s \| x^j(-\frac{1}{x}\frac{d}{dx})^j f \|_{L_2} \label{7.14}
\end{eqnarray}
задают эквивалентные нормировки в пространстве Соболева при целых $s \in \mathbb{Z}$. Кроме того, каждое слагаемое в \eqref{7.13} тождественно равно соответствующему слагаемому в \eqref{7.14} с тем же индексом $j$.

И. А. Киприянов ввёл в \cite{Kip2} шкалу пространств, которые оказали существенное влияние на теорию уравнений в частных производных с оператором Бесселя по одной или нескольким переменным. Эти пространства можно определить следующим образом. Рассмотрим подмножество чётных функций на $\mathbb{D}(0, \infty)$, у которых все производные нечётного порядка равны нулю при $x=0$. Обозначим это множество $\mathbb{D}_c (0, \infty)$ и введём на нём норму
\begin{equation}\label{7.15}{\|f\|_{\widetilde{W}_{2, k}^s} = \|f\|_{L_{2, k}}+\|B_k^{\frac{s}{2}}\|_{L_{2, k}}}
\end{equation}
где $s$ --- чётное натуральное число,  $B^{s/2}_k$ --- итерация оператора Бесселя. Пространство И.А.~Киприянова при чётных $s$ определяется как замыкание $D_c (0, \infty)$ по норме \eqref{7.15}. Известно, что эквивалентная \eqref{7.15} норма может быть задана по формуле \cite{Kip2}
\begin{equation}\label{7.16}{\|f\|_{\widetilde{W}_{2, k}^s} = \|f\|_{L_{2, k}}+\|x^s(-\frac{1}{x}\frac{d}{dx})^s f\|_{L_{2, k}}}
\end{equation}
Это позволяет доопределить норму в $\widetilde{W}_{2, \, k}^s$ для всех $s$. Отметим, что по существу этот подход совпадает с одним из принятых в \cite{Kip2}, другой подход основан на использовании преобразования Фурье--Бесселя. Далее будем считать, что $\widetilde{W}_{2, k}^s$ нормируется по формуле \eqref{7.16}.

Введём весовую соболевскую норму
\begin{equation}\label{7.17}{\|f\|_{W_{2, k}^s} = \|f\|_{L_{2, k}}+\|D_-^s f\|_{L_{2, k}}}
\end{equation}
и пространство $W_{2, \, k}^s$ как замыкание $\mathbb{D}_c (0, \infty)$ по этой норме.

\theor\label{tvloz1}
а) Пусть $k \neq -n, ~ n \in \mathbb{N}$. Тогда пространство  $\widetilde{W}_{2, \, k}^s$ непрерывно вложено в $W_{2, \, k}^s$, причём существует постоянная $A_3>0$ такая, что
\begin{equation}\label{7.18}{\|f\|_{W_{2, k}^s}\leq A_3 \|f\|_{\widetilde{W}_{2, k}^s},}
\end{equation}
б) Пусть $k+s \neq -2m_1-1, ~ k-s \neq -2m_2-2, ~ m_1 \in \mathbb{N}_0, ~ m_2 \in \mathbb{N}_0$. Тогда справедливо обратное вложение $W_{2, \, k}^s$ в $\widetilde{W}_{2, \, k}^s$, причём существует постоянная $A_4>0$, такая, что
\begin{equation}\label{7.19}{\|f\|_{\widetilde{W}_{2, k}^s}\leq A_4 \|f\|_{W_{2, k}^s}.}
\end{equation}
в) Если указанные условия не выполняются, то соответствующие вложения не имеют места.

Следствие 1.  Пусть выполнены условия: $k \neq -n, ~ n \in \mathbb{N}$; $k+s \neq -2m_1-1,  ~ m_1 \in \mathbb{N}_0; ~ k-s \neq -2m_2-2, ~ m_2 \in \mathbb{N}_0$. Тогда пространства И.А. Киприянова можно определить как замыкание $D_c (0, \infty)$ по весовой соболевской норме \eqref{7.17}.

Следствие 2.  Точные значения постоянных в неравенствах вложения \eqref{7.18}--\eqref{7.19} есть
$$
A_3 = \max (1, \|{_1S_-^{s-1}} \| _ {L_{2, k}}), ~ A_4=\max(1, \|{_1P_-^{s-1}}\|_{L_{2, k}}).
$$

Очевидно, что приведённая теорема и следствия из неё вытекают из приведённых выше результатов для операторов Бушмана--Эрдейи. Отметим, что  нормы операторов Бушмана--Эрдейи нулевого порядка гладкости в $L_{2, \, k}$  дают значения точных постоянных в неравенствах вложения \eqref{7.18}--\eqref{7.19}. Оценки норм операторов Бушмана--Эрдейи в банаховых пространствах $L_{p, \alpha}$ позволяют рассматривать вложения для соответствующих функциональных пространств.

Таким образом в этом пункте с помощью ОП Бушмана--Эрдейи нулевого порядка гладкости дан положительный ответ на вопрос, который давно обсуждался в устном "фольклоре"\ --- \textit{пространства И.А.Киприянова изоморфны весовым пространствам С.Л.Соболева}. Разумеется, мы рассмотрели самый простой случай, результаты можно обобщать на другие виды нормировок, многомерный случай, замену неограниченных областей на ограниченные, что ещё предстоит рассмотреть в будущем, но это не изменит принципиально основного вывода. Сказанное ни в коем случае не умаляет ни существенного значения, ни необходимости использования  пространств И.А.\,Киприянова для подходящего круга задач теории функций и дифференциальных уравнений с частными производными.
Принципиальная важность пространств И.А.\,Киприянова для теории уравнений в частных производных различных типов с операторами Бесселя отражает общий методологический подход, который автор услышал в виде красивого афоризма на пленарной лекции чл.--корр. РАН Л.Д.\,Кудрявцева:
\begin{center}
 "\textit{КАЖДОЕ УРАВНЕНИЕ ДОЛЖНО ИЗУЧАТЬСЯ В СВОЁМ СОБСТВЕННОМ ПРОСТРАНСТВЕ!}"
\end{center}

Доказанные в этом пункте вложения могут быть использованы для прямого переноса известных оценок для решений $B$--эллиптических уравнений в пространствах  И.А.\,Киприянова (см., например, \cite{Kip1},\cite{Kip2}) на оценки в  весовых пространствах С.Л.Соболева, это непосредственное применение приведённых в статье условий ограниченности и сплетающих свойств операторов преобразования Бушмана--Эрдейи.

\bigskip

\subsection{Представление решений дифференциальных уравнений в частных производных
с операторами Бесселя.}

Построенные операторы преобразования позволяют выписывать явные формулы, выражающие решения уравнений в частных производных с операторами Бесселя через невозмущённые уравнения. Примером служит
 $B$--эллиптическое уравнение  с операторами Бесселя по каждой переменной вида
\begin{equation}
\sum_{k=1}^{n}B_{\nu,x_k}u(x_1,\dots, x_n)=f,
\end{equation}
аналогичные  $B$--гиперболические и $B$--параболические уравнения. Эта идея ранее осуществлялась с использованием операторов преобразования Сонина--Пуассона--Дельсарта, см. \cite{Car1}--\cite{Car3}, \cite{CarSho}, \cite{Kip1}. Новые типы операторов преобразования позволяют получить новые классы подобных формул соответствия.

\bigskip

\subsection{Задача Коши для уравнения Эйлера--Пуассона--Дарбу (ЭПД).}

Рассмотрим уравнение ЭПД в полупространстве
$$
B_{\alpha,\, t} u(t,x)= \frac{ \partial^2 \, u }{\partial t^2} + \frac{2 \alpha+1}{t} \frac{\partial u}{\partial t}=\Delta_x u+F(t, x),
$$
где $t>0,~x \in \mathbb{R}^n$. Дадим нестрогое описание процедуры, позволяющей получать различные постановки начальных условий при $t=0$ единым методом. Образуем по формулам \eqref{75} операторы преобразования $X_{\alpha, \, t}$ и $Y_{\alpha, \, t}$. Предположим, что существуют выражения $X_{\alpha, \, t} u=v(t,x)$, $X_{\alpha, \, t} F=G(t,x)$. Пусть обычная (несингулярная) задача Коши
\begin{equation}\label{7.28}{\frac{ \partial^2 \, v }{\partial t^2} =\Delta_x v+G,~ v|_{t=0}=\varphi (x),~ v'_t|_{t=0}=\psi (x)}
\end{equation}
корректно разрешима в полупространстве. Тогда в предположении, что $Y_{\alpha, \, t}=X^{-1}_{\alpha, \, t}$ получаем следующие начальные условия для уравнения ЭПД:
\begin{equation}\label{7.29}{X_{\alpha} u|_{t=0}=a(x),~(X_{\alpha} u)'|_{t=0}=b(x).}
\end{equation}
При этом различному выбору операторов преобразования  $X_{\alpha, t}$ (операторы Сонина--Пуассона--Дельсарта, Бушмана--Эрдейи первого, второго, третьего родов, Бушмана--Эрдейи нулевого порядка гладкости, унитарные операторы преобразования Сонина--Катрахова и Пуассона--Катрахова, обобщенные операторы Бушмана--Эрдейи) будут соответствовать различные начальные условия. Следуя изложенной методике в каждом конкретном случае их можно привести к более простым аналитическим формулам, см. \cite{Sit3}.

Данная схема обобщается на дифференциальные уравнения с большим числом переменных, по которым могут действовать операторы Бесселя с различными параметрами, а также уравнения других типов.
Применение операторов преобразований позволяет сводить сингулярные (или иначе вырождающиеся) уравнения с операторами Бесселя по одной или нескольким переменным (уравнения ЭПД, сингулярное уравнение теплопроводности, $B$---эллиптические уравнения по определению И.А.\,Киприянова, уравнения обобщённой осесимметрической теории потенциала---теории $GASPT (Generalized\  Axially\  Symmetric\  Potential\  Theory)$---А.\,Вайнстейна и другие) к несингулярным. При этом априорные оценки для сингулярного случая получаются как следствия соответствующих априорных оценок для регулярных уравнений, если только удалось оценить сами операторы преобразования  в нужных функциональных пространствах. Значительное число подобных оценок было приведено выше.

Из работ этого направления отметим  монографию А.\,Псху (см. \cite{Psh}), в которой по существу применяется указанная выше схема с операторами преобразованиями для решения задачи Коши для уравнения с дробными производными, при этом существенно используется преобразование Станковича; а также применение операторов преобразования Бушмана--Эрдейи в работах А.В.\,Глушака, см. \cite{Glu}.

Отметим, что операторы Бушмана--Эрдейи и возникли впервые в теории уравнения Эйлера--Пуассона--Дарбу. Приведём соответствующий результат, имеющий в том числе исторический интерес.

\bigskip

Лемма Копсона.\\
Рассмотрим дифференциальное уравнение в частных производных с двумя переменными:

$$
\frac{\partial^2 u(x,y)}{\partial x^2}+\frac{2\alpha}{x}\frac{\partial u(x,y)}{\partial x}=
\frac{\partial^2 u(x,y)}{\partial y^2}+\frac{2\beta}{y}\frac{\partial u(x,y)}{\partial y}
$$
(обобщённое уравнение Эйлера--Пуассона--Дарбу или В--гиперболическое уравнение по терминологии И.А.\,Киприянова) в открытой четверти плоскости $x>0,\  y>0$ при положительных параметрах $\beta>\alpha>0$ с краевыми условиями на осях координат (характеристиках)
$$
u(x,0)=f(x), u(0,y)=g(y), f(0)=g(0).
$$
Предполагается, что решение u(x,y) является непрерывно дифференцируемым в замкнутом первом квадранте, имеет непрерывные вторые производные в открытом квадранте, граничные функции $f(x), g(y)$ являются непрерывно дифференцируемыми.

Тогда, если решение поставленной задачи существует, то для него выполняются соотношения:

\begin{equation}
\label{Cop1}
\frac{\partial u}{\partial y}=0, y=0,  \frac{\partial u}{\partial x}=0, x=0,
\end{equation}
\begin{equation}
\label{Cop2}
2^\beta \Gamma(\beta+\frac{1}{2})\int_0^1 f(xt)t^{\alpha+\beta+1}
\lr{1-t^2}^{\frac{\beta -1}{2}}P_{-\alpha}^{1-\beta}{t}\,dt=
\end{equation}
\begin{equation*}
=2^\alpha \Gamma(\alpha+\frac{1}{2})\int_0^1 g(xt)t^{\alpha+\beta+1}
\lr{1-t^2}^{\frac{\alpha -1}{2}}P_{-\beta}^{1-\alpha}{t}\,dt,
\end{equation*}
$$
\Downarrow
$$
\begin{equation}
\label{Cop3}
g(y)=\frac{2\Gamma(\beta+\frac{1}{2})}{\Gamma(\alpha+\frac{1}{2})
\Gamma(\beta-\alpha)}y^{1-2\beta}
\int_0^y x^{2\alpha-1}f(x)
\lr{y^2-x^2}^{\beta-\alpha-1}x \,dx,
\end{equation}
где  $P_\nu^\mu(z)$---функция Лежандра первого рода \cite{Sit1}.

Таким образом, содержание леммы Копсона сводится к тому, что начальные данные на характеристиках нельзя задавать произвольно, они должны быть связаны операторами Бушмана--Эрдейи первого рода. Более подробное обсуждение этой леммы и соответствующие ссылки см. в \cite{Sit1}.

\bigskip

\subsection{Применения к операторам обобщённого сдвига}.

Данный класс операторов введён и подробно изучен в работах Б.М.\,Левитана \cite{Lev1}--\cite{Lev2}.
Он имеет многочисленные применения в теории операторов с частными производными, в том числе с операторами Бесселя \cite{Lev3}, позволяя в частности переносить особенность в уравнениях из начала координат в любую точку. Операторы обобщённого сдвига по явным формулам выражаются через операторы преобразования  \cite{Lev3}. Поэтому новые классы операторов преобразования позволяют построить и изучать новые классы
операторов обобщённого сдвига.

\bigskip

\subsection{Применения к операторам Дункла. }

В последнее время значительное развитие получила теория операторов Дункла. Это в существенном дифференциально--разностные операторы, содержащие линейные комбинации обычных производных и конечных разностей. В высших размерностях операторы Дункла связаны с алгебрами Ли и группами отражений и симметрий.
Для этого класса операторов значительное развитие получила теория операторов преобразования, как классических, так в отдельных работах и Бушмана--Эрдейи, см., например, \cite{Sig} и ссылки в статьях этого сборника.

\bigskip

\subsection{Применения операторов преобразования Бушмана--Эрдейи в теории преобразования Радона. }

Известно, что в силу результатов Людвига \cite{Lud}  преобразование Радона    при описании через сферические гармоники действует на каждой гармонике по радиальной переменной в нашей терминологии как некоторый оператор Бушмана--Эрдейи первого рода. Приведём точный результат.

\theor Теорема Людвига \ (\cite{Lud},\cite{Hel}). Пусть задано разложение функции в $\mathbb{R}^n$ по сферическим гармоникам вида
\begin{equation}
f(x)=\sum_{k,l}f_{k,l}(r) Y_{k,l}(\theta).
\end{equation}
Тогда преобразование Радона этой функции также может быть разложено в ряд по сферическим гармоникам по формулам
\begin{equation}
Rf(x)=g(r,\theta)=\sum_{k,l}g_{k,l}(r) Y_{k,l}(\theta),
\end{equation}
\begin{equation}
\label{lu1}
g_{k,l}(r)=с(n)\int_r^\infty \lr{1-\frac{s^2}{r^2}}^{\frac{n-3}{2}} C_l^{\frac{n-2}{2}}\lr{\frac{s}{r}}
f_{k,l}(r) r^{n-2}\,ds,
\end{equation}
где $с(n)$---некоторая известная постоянная, $C_l^{\frac{n-2}{2}}\lr{\frac{s}{r}}$---функция Гегенбауэра \cite{BE1}. Также справедлива аналогичная обратная формула, выражающая величины $f_{k,l}(r)$ через $g_{k,l}(r)$.

Известно, что функция Гегенбауэра выражается через функцию Лежандра \cite{BE1}. Таким образом, формула Людвига (\ref{lu1}) показывает, что результат действия преобразования Радона по каждой гармонике с точностью до степенных и постоянных множителей---это операторы преобразования Бушмана--Эрдейи первого рода.

Отметим, что именно эта формула в двумерном случае была использована А.\,Кормаком для расчёта первого томографа, за что впоследствии он был удостоен Нобелевской премии.

Частными случаями формулы Людвига, полученной в 1966 году, являются явные формулы, описывающие действие преобразования Радона по любой сферической гармонике, в частности, на чисто радиальных функциях. Эти формулы переоткрываются энтузиастами до сих пор.

На этом пути следствием полученных выше результатов являются интегральные представления, оценки норм в функциональных пространствах, формулы обращения для преобразования Радона. Результаты формулируются в терминах операторов Сонина--Пуассона--Дельсарта, Эрдейи--Кобера, Бушмана--Эрдейи.

Например, становится понятным, что  по существу многие формулы обращения преобразования Радона совпадают с различными вариантами формул для обращения операторов Бушмана--Эрдейи первого рода.

Указанный круг вопросов также подробно изложен в монографии \cite{Dea}.

\bigskip

\subsection{ Применения операторов преобразования Бушмана--Эрдейи к построению обобщённых сферических гармоник.}

Ещё с 1950--х годов специалистам было известно, что никакой новой теории для построения полиномиальных решений В--эллиптических уравнений не требуется. Это следует из простого факта, что уже операторы преобразования  Сонина--Пуассона--Дельсарта переводят степень в степень. Следовательно, они переводят В--гармонические полиномы в гармонические и наоборот по явным формулам. Поэтому и обобщённые сферические гармоники строятся по явным формулам из обычных, так как они являются сужениями соответствующих полиномов на единичную сферу. Данный подход подробно изложен, например, в работах Б.Рубина  \cite{Rub1}--\cite{Rub2} (см. также ссылки в этих работах). Операторы Бушмана--Эрдейи добавляют новую степень свободы ко всем этим построениям.

\bigskip

\subsection{Применения операторов преобразования Бушмана--Эрдейи к построению унитарных обобщений операторов Харди.}

Унитарность сдвинутых на единичный классических операторов Харди (\ref{3.25}) установлена выше, см. следствие 4 из теоремы 7, как уже отмечалось---это известный результат, приведённый в \cite{KuPe}. Подстановка других натуральных значений параметра приводит к новому бесконечному семейству интегральных операторов простого вида, унитарных в $L_2(0,\infty)$.

\theor Следующие операторы образуют пары взаимно обратных унитарных операторов в  $L_2(0,\infty)$:
\begin{eqnarray*}
\label{84}
U_3f= f+\int_0^x f(y)\,\frac{dy}{y},\  U_4f= f+\frac{1}{x}\int_x^\infty f(y)\,dy,\\\nonumber
U_5f= f+3x\int_0^x f(y)\,\frac{dy}{y^2},\  U_6f= f-\frac{3}{x^2}\int_0^x y f(y)\,dy,\\\nonumber
U_7f=f+\frac{3}{x^2}\int_x^\infty y f(y)\,dy,\  U_8f=f-3x \int_x^\infty f(y)\frac{dy}{y^2},\\\nonumber
U_9f=f+\frac{1}{2}\int_0^x \left(\frac{15x^2}{y^3}-\frac{3}{y}\right)f(y)\,dy,\\\nonumber
U_{10}f=f+\frac{1}{2}\int_x^\infty \left(\frac{15y^2}{x^3}-\frac{3}{x}\right)f(y)\,dy.\\\nonumber
\end{eqnarray*}

\bigskip

\subsection{Интегральные операторы с более общими специальными функциями в ядрах.}

Рассмотрим оператор ${{}_1S_{0+}^{\nu}}$. Он имеет вид
\begin{equation}\label{7.30}{{{_1}S_{0+}^{\nu}}=\frac{d}{dx}\int\limits_0^x K\lr{\frac{x}{y}} f(y) \, dy,}
\end{equation}
где ядро $K$ выражается по формуле $K(z)=P_{\nu}(z)$. Простейшие свойства специальных функций позволяют показать, что ${{_1}S_{0+}^{\nu}}$ можно рассматривать как частный случай оператора вида \eqref{7.30} с функцией Гегенбауэра в ядре
\begin{equation}\label{7.31}{K(z)=\frac{\Gamma(\alpha+1)\ \Gamma(2 \beta)}{2^{p-\frac{1}{2}}\Gamma(\alpha+2\beta)\ \Gamma(\beta+ \frac{1}{2})}(z^{\alpha}-1)^{\beta-\frac{1}{2}}C^{\beta}_{\alpha}(z)}
\end{equation}
при значениях параметров $\alpha=\nu, ~ \beta = \frac{1}{2}$ или с функцией Якоби в ядре
\begin{equation}\label{7.32}{K(z)=\frac{\Gamma(\alpha+1)}{2^{\rho}\Gamma(\alpha+\rho+1)}(z-1)^{\rho}(z+1)^{\sigma}P^{(\rho, \sigma)}_{\alpha}(z)}
\end{equation}
при значениях параметров $\alpha=\nu, ~ \rho=\sigma = 0$. Более общим являются операторы с гипергеометрической функцией Гаусса ${{_2}F_1}$ или с $G$ --- функцией  Майера в ядрах. Перейдём к их определению.

Операторы с функцией Гаусса ${{_2}F_1}(a, b; c; z)$ изучались в большом числе работ, см. библиографию в \cite{KK}.   Подобные обобщения другого рода, в которых рассматриваются операторы с интегрированием по всей полуоси и ядра  выражаются через обобщённые функции Лежандра, изучались в \cite{ViFe}.
Также выделим монографию А.А.\,Килбаса и М.\,Сайго \cite{KiSa}, в которой подробно изложены интегральные преобразования с $H$--функцией Фокса в ядрах.
Для исследования таких операторов оказываются  полезными различные неравенства для гипергеометрических функций, например, полученные в \cite{SiKa1}--\cite{SiKa2}.

Введём ещё один класс подобных операторов, обобщающих операторы Бушмана--Эрдейи \eqref{71}--\eqref{72}.

Definition 6.{
 Определим операторы Гаусса--Бушмана--Эрдейи по  следующим формулам:
\begin{equation}\label{7.33}{{_1F_{0+}}(a, b, c)[f]=\frac{1}{2^{c-1}\Gamma(c)}.}
\end{equation}
$$
\int\limits_0^x\lr{\frac{x}{y}-1}^{c-1}\lr{\frac{x}{y}+1}^{a+b-c} {{_2}F_1}\lr{^{a,b}_c| \frac{1}{2}-\frac{1}{2}\frac{x}{y}} f(y) \, dy,
$$
\begin{equation}\label{7.34}{{{_2}F_{0+}}(a, b, c)[f]=\frac{1}{2^{c-1}\Gamma(c)}.}
\end{equation}
$$
\int\limits_0^x\lr{\frac{y}{x}-1}^{c-1}\lr{\frac{y}{x}+1}^{a+b-c} {{_2}F_1}\lr{^{a,b}_c| \frac{1}{2}-\frac{1}{2}\frac{y}{x}} f(y) \, dy,
$$
\begin{equation}\label{7.35}{{_1F_{-}}(a, b, c)[f]=\frac{1}{2^{c-1}\Gamma(c)}.}
\end{equation}
$$
\int\limits_0^x\lr{\frac{y}{x}-1}^{c-1}\lr{\frac{y}{x}+1}^{a+b-c} {{_2}F_1}\lr{^{a,b}_c| \frac{1}{2}-\frac{1}{2}\frac{y}{x}} f(y) \, dy,
$$
\begin{equation}\label{7.36}{{{_2}F_{-}}(a, b, c)[f]=\frac{1}{2^{c-1}\Gamma(c)}.}
\end{equation}
$$
\int\limits_0^x\lr{\frac{x}{y}-1}^{c-1}\lr{\frac{x}{y}+1}^{a+b-c} {{_2}F_1}\lr{^{a,b}_c| \frac{1}{2}-\frac{1}{2}\frac{x}{y}} f(y) \, dy,
$$
\begin{eqnarray}
& {{_3}F_{0+}}[f]=\frac{d}{dx} {{_1}F_{0+}}[f], & {{_4}F_{0+}}[f]= {{_2}F_{0+}} \frac{d}{dx} [f], \label{7.37} \\
& {{_3}F_{-}}[f]={{_1}F_{-}}(-\frac{d}{dx} ) [f], & {{_4}F_{-}}[f]= (- \frac{d}{dx} ) {{_2}F_{-}}[f]. \label{7.38}
\end{eqnarray}
}

Символ ${{_2}F_1}$  в определениях \eqref{7.34} и \eqref{7.36} означает гипергеометрическую функцию Гаусса на естественной области определения, а в определениях \eqref{7.33} и \eqref{7.35}
обозначает главную ветвь аналитического продолжения этой функции.

Операторы \eqref{7.33}--\eqref{7.36} обобщают операторы Бушмана--Эрдейи \eqref{71}--\eqref{72} соответственно. Они сводятся к последним при выборе параметров $a=-(\nu+\mu),~ b= 1+ \nu - \mu,~ c=1-\mu$. На операторы \eqref{7.33}--\eqref{7.36} с соответствующими изменениями переносятся  все полученные выше результаты. В частности они факторизуются через более простые операторы \eqref{7.37}--\eqref{7.38} при специальном выборе параметров.

Операторы \eqref{7.37}--\eqref{7.38} обобщают операторы \eqref{73}--\eqref{733}. Для них справедлива

\theor Операторы \eqref{7.37}--\eqref{7.38} могут быть расширены до изометричных в $L_2 (0, \infty)$ тогда и только тогда, когда они совпадают с операторами Бушмана--Эрдейи нулевого порядка гладкости I рода \eqref{73}--\eqref{733} соответственно при целых значениях параметра $\nu=\frac{1}{2}(b-a-1)$.

Эта теорема выделяет операторы Бушмана--Эрдейи нулевого порядка гладкости среди их возможных обобщений по крайней мере вида \eqref{7.33}--\eqref{7.38}. Сами операторы
\eqref{7.33}--\eqref{7.36} интересны как новый класс преобразований, обобщающих операторы дробного интегродифференцирования. Аналогичные обобщения можно проделать и для операторов \eqref{6.1}--\eqref{6.2}, \eqref{6.6}, \eqref{6.16}--\eqref{6.17}.

Более общими являются операторы с $G$--функцией Майера в ядре. Например, один из таких операторов имеет вид

\begin{eqnarray}
& & {{_1}G_{0+}}(\alpha, \beta, \delta, \gamma)[f]= \frac{2^{\delta}}{\Gamma(1-\alpha)\Gamma(1-\beta)}\cdot \label{7.39} \\
& & \int\limits_0^x (\frac{x}{y}-1)^{-\delta}(\frac{x}{y}+1)^{1+\delta-\alpha-\beta} G_{2~2}^{1~2} \lr{\frac{x}{2y}-\frac{1}{2}|_{\gamma, \, \delta}^{\alpha, \, \beta}}  f(y)\, dy. \nonumber
\end{eqnarray}
Остальные получаются при изменении промежутка интегрирования и значений аргумента $G$ -- функции. При значениях параметров $\alpha=1-a,~ \beta=1-b, \delta=1-c, \gamma=0$ \eqref{7.39} сводится к \eqref{7.33}, а при значениях $\alpha=1+\nu,~ \beta=-\nu, \delta= \gamma=0$  \eqref{7.39} сводится к оператору Бушмана--Эрдейи I рода нулевого порядка гладкости ${{_1}S_{0+}^{\nu}}$.

Дальнейшие обобщения возможны в терминах функций Райта или Фокса. На этом пути после операторов преобразования Бушмана--Эрдейи и Гаусса--Бушмана--Эрдейи вводятся операторы преобразования Майера--Бушмана--Эрдейи,
Райта--Бушмана--Эрдейи и Фокса--Бушмана--Эрдейи различных родов. Отметим, что вводимые конструкции связаны с операторами преобразования для гипербесселевых дифференциальных операторов типа Сонина--Димовски и Пуассона--Димовски \cite{Dim},  \cite{Kir},   а также обобщёнными интегродифференциальными операторами дробного порядка, которые были введены В.\,Киряковой \cite{Kir}.

\bigskip

\subsection{Применение операторов преобразования Бушмана--Эрдейи в работах В.В.\,Катрахова.}

В.В.\,Катраховым был предложен новый подход к постановке краевых задач для эллиптических уравнений с особенностями. Например, для уравнения Пуассона им рассматривалась задача в области, содержащей начало координат, в которой решения могут иметь особенности произвольного роста. В точке начала координат им было предложено новое нелокальное краевое условие типа свертки, которое мы назовём К---следом. В определении классов для решений, которые обобщают пространства С.Л.\,Соболева на случай функций с существенными особенностями, фундаментальную роль играют различные ОП. Основные результаты состоят в доказательстве корректной разрешимости поставленных задач во введённых пространствах. И постановки задач с К---следом, и использованные нормы напрямую используют операторы преобразования Бушмана--Эрдейи, см. \cite{Kat1}--\cite{Kat2}.

Кроме того, в совместных работах  И.А.\,Киприянова и В.В.\,Катрахова на основе операторов преобразования вводятся и изучаются новые классы ПДО \cite{Kat3}. Эти результаты изложены в качестве отдельного параграфа в монографии Р.\,Кэрролла \cite{Car2}.

\bigskip

\begin{center}
ЛИТЕРАТУРА.
\end{center}

\addcontentsline{toc}{section}{ЛИТЕРАТУРА \hspace{99mm} 61}

\selectlanguage{russian}

\tableofcontents

\end{document}